\RequirePackage{ifpdf}
\ifpdf % We are running pdfTeX in pdf mode
\documentclass[pdftex]{sigma}
\else
\documentclass{sigma}
\fi

\numberwithin{equation}{section}

\newcommand{\Zint}{\mathbb {Z}}

\newcommand{\Rea}{\mathbb {R}}
\newcommand{\Cplx}{\mathbb {C}}

\begin{document}

\allowdisplaybreaks

\renewcommand{\thefootnote}{$\star$}

\renewcommand{\PaperNumber}{040}

\FirstPageHeading

\ShortArticleName{Middle Convolution and Heun's Equation}

\ArticleName{Middle Convolution and Heun's Equation\footnote{This paper is a contribution to the Proceedings of the Workshop ``Elliptic Integrable Systems, Isomonodromy Problems, and Hypergeometric Functions'' (July 21--25, 2008, MPIM, Bonn, Germany). The full collection
is available at
\href{http://www.emis.de/journals/SIGMA/Elliptic-Integrable-Systems.html}{http://www.emis.de/journals/SIGMA/Elliptic-Integrable-Systems.html}}}

\Author{Kouichi TAKEMURA}

\AuthorNameForHeading{K. Takemura}

\Address{Department of Mathematical Sciences, Yokohama City University,\\ 22-2 Seto, Kanazawa-ku, Yokohama 236-0027, Japan}

\Email{\href{mailto:takemura@yokohama-cu.ac.jp}{takemura@yokohama-cu.ac.jp}}

\ArticleDates{Received November 26, 2008, in f\/inal form March 25,
2009; Published online April 03, 2009}

\Abstract{Heun's equation naturally appears as special cases of Fuchsian system of dif\/ferential equations of rank two with four singularities by introducing the space of initial conditions of the sixth Painlev\'e equation.
Middle convolutions of the Fuchsian system are related with an integral transformation of Heun's equation.}


\Classification{34M35; 33E10; 34M55}

\section{Introduction}

Heun's equation is a standard form of a second-order Fuchsian dif\/ferential equation with four singularities, and it is given by
\begin{gather}
\frac{d^2y}{dz^2} + \left( \frac{\gamma}{z}+\frac{\delta }{z-1}+\frac{\epsilon}{z-t}\right) \frac{dy}{dz} +\frac{\alpha \beta z -q}{z(z-1)(z-t)} y=0,
\label{eq:Heun}
\end{gather}
with the condition
\begin{gather*}
\gamma +\delta +\epsilon =\alpha +\beta +1.
%\label{Heuncond}
\end{gather*}
The parameter $q$ is called an accessory parameter.
Although the local monodromy (local expo\-nent) is independent of $q$, the global monodromy (e.g.\ the monodromy on the cycle enclosing two singularities) depends on $q$.
Some properties of Heun's equation are written in the books~\cite{Ron,SL},
but an important feature related with the theory of f\/inite-gap potential for the case $\gamma, \delta, \epsilon , \alpha -\beta \in \Zint +\frac{1}{2}$ (see \cite{GW,Smi,Tak1,Tak2,Tak3,Tak4,Tak5,TV} etc.), which leads to an algorithm to calculate the global monodromy explicitly for all $q$, is not written in these books.

The sixth Painlev\'e equation is a non-linear ordinary dif\/ferential equation written as
\begin{gather}
  \frac{d^2\lambda }{dt^2} =  \frac{1}{2} \left( \frac{1}{\lambda }+\frac{1}{\lambda -1}+\frac{1}{\lambda -t} \right) \left( \frac{d\lambda }{dt} \right) ^2 -\left( \frac {1}{t} +\frac {1}{t-1} +\frac {1}{\lambda -t} \right)\frac{d\lambda }{dt} \nonumber \\
\phantom{\frac{d^2\lambda }{dt^2} =}{}  +\frac{\lambda (\lambda -1)(\lambda -t)}{t^2(t-1)^2}\left\{ \frac{(1-\theta _{\infty})^2}{2} -\frac{\theta _{0}^2}{2}\frac{t}{\lambda ^2} +\frac{\theta  _{1}^2}{2}\frac{(t-1)}{(\lambda -1)^2} +\frac{(1-\theta _{t}^2)}{2}\frac{t(t-1)}{(\lambda -t)^2} \right\} . \label{eq:P6eqn}
\end{gather}

A remarkable property of this dif\/ferential equation is that the solutions do not have movable singularities other than poles.
It is known that the sixth Painlev\'e equation is obtained by monodromy preserving deformation of Fuchsian system of dif\/ferential equations,
\begin{gather*}
\frac{d}{dz}\left( \! \!
\begin{array}{c}
y_{1} \\
y_{2}
\end{array}
\! \! \right)= \left( \frac{A_0}{z}+\frac{A_1}{z-1}+\frac{A_t}{z-t} \right)
\left( \! \!
\begin{array}{c}
y_{1} \\
y_{2}
\end{array}
\! \! \right), \qquad A_0, A_1 , A_t \in \Cplx ^{2\times 2}.
%\label{eq:dy12dzAzy12}
\end{gather*}
See Section~\ref{sec:linP6} for expressions of the elements of the matrices $A_0$, $A_1$, $A_t$.
By eliminating $y_2$ we have second-order dif\/ferential equation for $y_1$, which have an additional apparent singularity $z=\lambda $ other than $\{ 0,1,t,\infty \}$ for generic cases, and the point $\lambda $ corresponds to the variable of the sixth Painlev\'e equation.
For details of monodromy preserving deformation, see~\cite{IKSY}.
In this paper we investigate the condition that the second-order dif\/ferential equation for $y_1$ is written as Heun's equation.
To get a preferable answer, we introduce the space of initial conditions for the sixth Painlev\'e equation which was discovered by Okamoto \cite{Oka} to construct a suitable def\/ining variety for the set of solutions to the (sixth) Painlev\'e equation.

For Fuchsian systems of dif\/ferential equations and local systems on a punctured Riemann sphere, Dettweiler and Reiter \cite{DR1,DR2} gave an algebraic analogue of Katz' middle convolution functor \cite{Katz}.
Filipuk~\cite{Fil} applied them for the Fuchsian systems with four singularities, obtained an explicit relationship with the symmetry of the sixth Painlev\'e equation, and the author~\cite{TakI} calculated the corresponding integral transformation for the Fuchsian systems with four singularities.
The middle convolution is labeled by a parameter $\nu $, and we have two values which leads to non-trivial transformation on $2\times 2$ Fuchsian system with four singularities (see Section~\ref{sec:MC}).
In this paper we consider the middle convolution which is a dif\/ferent value of the parameter $\nu$ from the one discussed in~\cite{Fil,TakI}.
We will also study the relationship between middle convolution and Heun's equation.
For special cases, the integral transformation raised by the middle convolution turns out to be a transformation on Heun's equation, and we investigate these cases.
Note that the description by the space of initial conditions for the sixth Painlev\'e equation is favorable.
The integral transformation of Heun's equation is applied for the study of novel solutions, which we will discuss in a separated publication.
If the parameter of the middle convolution is a negative integer, then the integral transformation changes to a successive dif\/ferential, and a transformation def\/ined by a dif\/ferential operator on Heun's equation was found in~\cite{Tak5} as a~generalized Darboux transformation (Crum--Darboux transformation).
Hence the integral transformation on Heun's equation can be regarded as a generalization of the generalized Darboux transformation, which is related with the conjectual duality by Khare and Sukhatme \cite{KS0}.

Special functions of the isomonodromy type including special solutions to the sixth Painlev\'e equation have been studied actively and they are related with various objects in mathematics and physics~\cite{Kit,VK}.
On the other hand, special functions of Fuchsian type including special solutions to Heun's equation are also interesting objects which are related with general relativity and so on.
This paper is devoted to an attempt to clarify both sides of viewpoints.

This paper is organized as follows:
In Section~\ref{sec:linP6}, we f\/ix notations for the Fuchsian system with four singularities.
In Section~\ref{sec:HeunSIC}, we def\/ine the space of initial conditions for the sixth Painlev\'e equation and observe that Heun's equation is obtained from the Fuchsian equation by restricting to certain lines in the space of initial conditions.
In Section~\ref{sec:MC}, we review results on the middle convolution and construct integral transformations.
In Section~\ref{sec:MCHeun}, we investigating relationship among the middle convolution, integral transformations of Heun's equation and the space of initial conditions.
In Section~\ref{sec:integer}, we consider the case that the parameter on the middle convolution is integer.
In the appendix, we describe topics which was put of\/f in the text.

\section{Fuchsian system of rank two with four singularities} \label{sec:linP6}

We consider a system of ordinary dif\/ferential equations,
\begin{gather}
\frac{dY}{dz}=A(z)Y, \!\qquad A(z)=\frac{A_0}{z}+\frac{A_1}{z-1}+\frac{A_t}{z-t} =
\left(\!\!
\begin{array}{ll}
a_{11}(z) & a_{12}(z) \\
a_{21}(z) & a_{22}(z)
\end{array}\!\!
\right) , \!\qquad Y= \left(\!\!
\begin{array}{l}
y_{1} \\
y_{2}
\end{array}\!\!
\right) ,\!\!\!
\label{eq:dYdzAzY1}
\end{gather}
where $t\neq 0,1$, $A_0$, $A_1$, $A_t$ are $2 \times 2$ matrix with constant elements.
Then equation~(\ref{eq:dYdzAzY1}) is Fuchsian, i.e., any singularities on the Riemann sphere $\Cplx \cup \{ \infty \} $ are regular, and it may have regular singularities at $z=0,1,t,\infty$ on the Riemann sphere $\Cplx \cup \{ \infty \}$.
Exponents of equation~(\ref{eq:dYdzAzY1}) at $z=0$ (resp. $z=1$, $z=t$, $z=\infty$) are described by eigenvalues of the matrix~$A_0$ (resp.~$A_1$, $A_t$, $-(A_0+A_1+A_t)$).
By the transformation $Y \rightarrow z^{n_0}(z-1)^{n_1}(z-t)^{n_t}Y $, the system of dif\/ferential equations (\ref{eq:dYdzAzY1}) is replaced as $A(z) \rightarrow A(z)+(n_0/z+n_1/(z-1)+n_2/(z-t))I$ ($I$: unit matrix),
and we can transform equation~(\ref{eq:dYdzAzY1}) to the one where one of the eigenvalues of $A_i$ is zero for $i\in \{ 0,1,t\}$ by putting $-n_i$ to be one of the eigenvalues of the original $A_i$.
If the exponents at $z=\infty $ are distinct, then we can normalize the matrix $-(A_0+A_1+A_t)$ to be diagonal by a suitable gauge transformation $Y \rightarrow GY$, $A(z) \rightarrow G A(z) G^{-1}$.
In this paper we assume that one of the eigenvalues of $A_i$ is zero for $i=0,1,t$ and the matrix $-(A_0 +A_1 +A_t) $ is diagonal, and we set
\begin{gather}
A_{\infty }= -(A_0 +A_1 +A_t) = \left(
\begin{array}{cc}
\kappa _1 & 0 \\
0 & \kappa _2
\end{array}
\right) . \label{def:Ainf}
\end{gather}
By eliminating $y_2$ in equation~(\ref{eq:dYdzAzY1}), we have a second-order linear dif\/ferential equation,
\begin{gather}
  \frac{d^2y_1}{dz^2} + p_1(z)  \frac{dy_1}{dz}+ p_2(z) y_1=0, \qquad p_1(z)= -a_{11}(z)-a_{22}(z)-\frac{\frac{d}{dz}a_{12}(z)}{a_{12}(z)}, \nonumber \\
  p_2 (z) = a_{11}(z)a_{22}(z)-a_{12}(z)a_{21}(z)- { {\frac{d}{dz}a_{11}(z)}}+\frac{a_{11}(z)\frac{d}{dz}a_{12}(z)}{a_{12}(z)} .  \label{eq:y1DE}
\end{gather}
Set
\begin{gather}
A_i =
\left(
\begin{array}{ll}
a_{11} ^{(i)} & a_{12}^{(i)} \\
a_{21} ^{(i)}& a_{22} ^{(i)}
\end{array}
\right) ,\qquad (i=0,1,t). \label{eq:Aiajk}
\end{gather}
It follows from equation~(\ref{def:Ainf}) that $a_{12}^{(0)}+ a_{12}^{(1)}+ a_{12}^{(t)}=0$, $a_{21}^{(0)}+ a_{21}^{(1)}+ a_{21}^{(t)}=0$. Hence $a_{12} (z)$ and~$a_{21} (z)$ are expressed as
\begin{gather*}
  a_{12} (z) =\frac{k_1 z +k_2 }{z(z-1)(z-t)}, \qquad a_{21} (z)=\frac{\tilde{k}_1 z +\tilde{k}_2 }{z(z-1)(z-t)},
%\label{eq:a12a21k1tk1}
\end{gather*}
and we have
\begin{gather*}
  a_{12}^{(0)}+ a_{12}^{(1)}+ a_{12}^{(t)}=0, \qquad (t+1)a_{12}^{(0)}+ t a_{12}^{(1)}+ a_{12}^{(t)}=-k_1 ,\qquad ta_{12}^{(0)}= k_2 , \\
  a_{21}^{(0)}+ a_{21}^{(1)}+ a_{21}^{(t)}=0, \qquad (t+1)a_{21}^{(0)}+ t a_{21}^{(1)}+ a_{21}^{(t)}=-\tilde{k}_1 ,\qquad ta_{21}^{(0)}= \tilde{k}_2 . \nonumber
\end{gather*}
If $k_1=k_2=0$, then $y_1$ satisf\/ies a f\/irst-order linear dif\/ferential equation, and it is integrated easily.
Hence we assume that $(k_1,k_2) \neq (0,0)$.
Then it is shown that two of $a_{12}^{(0)}$, $a_{12}^{(1)}$, $a_{12}^{(t)}$, $(t+1)a_{12}^{(0)}+ t a_{12}^{(1)}+ a_{12}^{(t)}$ cannot be zero.
We set $\lambda = -k_2/k_1$ ($k_1\neq 0$) and $\lambda = \infty$ ($k_1 = 0$).
The condition that none of $a_{12}^{(0)}$, $a_{12}^{(1)}$, $a_{12}^{(t)}$ nor $(t+1)a_{12}^{(0)}+ t a_{12}^{(1)}+ a_{12}^{(t)}$ is zero is equivalent to that $\lambda \neq 0,1,t, \infty$,
and the condition $a_{12}^{(0)} =0$ (resp.\  $a_{12}^{(1)}=0$, $a_{12}^{(t)}=0$, $(t+1)a_{12}^{(0)}+ t a_{12}^{(1)}+ a_{12}^{(t)}=0$)  is equivalent to $\lambda = 0$ (resp.\ $\lambda = 1$, $\lambda = t$, $\lambda = \infty$).

We consider the case $\lambda \neq 0,1,t, \infty$, i.e., the case $a_{12}^{(0)}\neq 0$, $a_{12}^{(1)}\neq 0$, $a_{12}^{(t)}\neq 0$, $(t+1)a_{12}^{(0)}+ t a_{12}^{(1)}+ a_{12}^{(t)}\neq 0 $.
Let $\theta _0$ (resp.\ $\theta _1$, $\theta _t$) and $0$ be the eigenvalues of $A_0$ (resp.\ $A_1$, $A_t$).
Then we can set $A_0$, $A_1$, $A_t$ as
\begin{gather}
 A_0= \left(
\begin{array}{ll}
u_0+\theta _0  & -w_0 \\
u_0(u_0+\theta _0)/ w_0  & -u_0
\end{array}
\right) , \qquad
A_1= \left(
\begin{array}{ll}
u_1+\theta _1  & -w_1 \\
u_1(u_1+\theta _1)/w_1 & -u_1
\end{array}
\right) , \nonumber \\
  A_t= \left(
\begin{array}{ll}
u_t+\theta _t  & -w_t \\
u_t(u_t+\theta _t)/w_t  & -u_t
\end{array}
\right) ,    \label{eq:A0A1AtP}
\end{gather}
by introducing variables $u_0$, $w_0$, $u_1$, $w_1$, $u_t$, $w_t$.
By taking trace of equation~(\ref{def:Ainf}), we have the relation $\theta _0 +\theta _1+\theta _t +\kappa _1+\kappa _2 =0$. We set $\theta _{\infty }=\kappa _1-\kappa _2$, then we have $\kappa _1= (\theta _{\infty } -\theta _0 -\theta _1 -\theta _t)/2$, $\kappa _2= -(\theta _{\infty } +\theta _0 +\theta _1 +\theta _t)/2$.

We determine $u_0$,  $u_1$,  $u_t$,  $w_0$,  $w_1$,  $w_t$ so as to satisfy equation~(\ref{def:Ainf}) and the following relations:
\begin{gather*}
  a_{12} (z) =-\frac{w_0}{z}-\frac{w_1}{z-1}-\frac{w_t}{z-t} =\frac{k( z -\lambda ) }{z(z-1)(z-t)}, \\
  a_{11}(\lambda)=\frac{u_0+\theta _0}{\lambda}+\frac{u_1+\theta _1}{\lambda-1}+\frac{u_t+\theta _t}{\lambda-t} =\mu , \nonumber
\end{gather*}
(see \cite{JM}). Namely, we solve the following equations for $u_0$,  $u_1$,  $u_t$,  $w_0$,  $w_1$,  $w_t$:
\begin{gather}
  -w_0-w_1-w_t=0, \qquad w_0(t+1)+w_1 t+w_t =k, \qquad -w_0 t=-k\lambda ,  \nonumber\\
  u_0(u_0+\theta _0)/ w_0 +u_1(u_1+\theta _1)/ w_1 +u_t(u_t+\theta _t)/ w_t =0, \nonumber \\
  u_0+\theta _0 + u_1+\theta _1 + u_t+\theta _t =-\kappa _1, \qquad  -u_0-u_1- u_t =-\kappa _2 , \nonumber \\
  (u_0+\theta _0)/\lambda+(u_1+\theta _1)/(\lambda-1)+(u_t+\theta _t)/(\lambda-t) =\mu .\label{eq:wuA01t}
\end{gather}
The linear equations for $w_0$, $w_1$, $w_t$ are solved as
\begin{gather}
  w_0 = \frac{k\lambda}{t}, \qquad  w_1= -\frac{k(\lambda-1)}{t-1}, \qquad w_t = \frac{k(\lambda-t)}{t(t-1)} . \label{eq:wgen}
\end{gather}
By the equations which are linear in $u_0$, $u_1$ and $u_t$, we can express $u_1 + \theta _1$ and $u_t +\theta _t$ as linear functions in $u_0$.
We substitute $u_1 + \theta _1$ and $u_t +\theta _t$ into a quadratic equation in $u_0$, $u_1$ and $u_t$.
Then the coef\/f\/icient of $u_0^2$ disappears, and $u_0$, $u_1$, $u_t$ are solved as
\begin{gather}
  u_0=-\theta _0 +\frac{\lambda}{t\theta _{\infty}} [ \lambda(\lambda-1)(\lambda-t)\mu^2 +\{ 2\kappa _1 (\lambda-1)(\lambda-t)-\theta _1(\lambda-t)  \nonumber\\
\phantom{u_0=}{}  -t\theta _t(\lambda-1) \} \mu +\kappa _1 \{\kappa _1(\lambda-t-1)-\theta _1-t\theta _t\} ] ,\nonumber \\
  u_1 =-\theta _1 -\frac{\lambda-1}{(t-1)\theta _{\infty}}  [ \lambda(\lambda-1)(\lambda-t)\mu^2 + \{ 2\kappa _1 (\lambda-1)(\lambda-t)+(\theta _{\infty}-\theta _1 )(\lambda-t) \nonumber\\
\phantom{u_1 =}{}  -t\theta _t (\lambda-1)\}\mu +\kappa _1 \{ \kappa _1 (\lambda-t+1)+\theta _0-(t-1)\theta _t\} ] ,\nonumber \\
  u_t= -\theta _t+\frac{\lambda-t}{t(t-1)\theta _{\infty}}  [ \lambda(\lambda-1)(\lambda-t)\mu^2 + \{ 2\kappa _1 (\lambda-1)(\lambda-t) -\theta _1 (\lambda-t) \nonumber\\
\phantom{u_t=}{}  +t(\theta _{\infty}-\theta _t )(\lambda-1)\} \mu +\kappa _1 \{ \kappa _1 (\lambda-t+1)+\theta _0+(t-1)(\theta _{\infty}-\theta _t)\} ] .\label{eq:ugen}
\end{gather}
We denote the Fuchsian system of dif\/ferential equations
\begin{gather}
\frac{dY}{dz}=\left( \frac{A_0}{z}+\frac{A_1}{z-1}+\frac{A_t}{z-t} \right) Y, \qquad Y= \left(
\begin{array}{l}
y_{1} \\
y_{2}
\end{array}
\right) ,
\label{eq:dYdzAzY}
\end{gather}
with equations~(\ref{eq:A0A1AtP}), (\ref{eq:wgen}), (\ref{eq:ugen}) by $D_Y(\theta _0, \theta _1, \theta _t, \theta _{\infty}; \lambda ,\mu ;k )$.
Then the second-order dif\/ferential equation (\ref{eq:y1DE}) is written as
\begin{gather}
  \frac{d^2y_1}{dz^2} + \left( \frac{1-\theta _0}{z}+\frac{1-\theta _1}{z-1}+\frac{1-\theta _t}{z-t}-\frac{1}{z-\lambda} \right) \frac{dy_1}{dz}  \nonumber\\
\phantom{\frac{d^2y_1}{dz^2}}{}  + \left( \frac{\kappa _1(\kappa _2 +1)}{z(z-1)}+\frac{\lambda (\lambda -1)\mu}{z(z-1)(z-\lambda)}-\frac{t (t -1)H}{z(z-1)(z-t)} \right) y_1=0, \nonumber \\
  H=\frac{1}{t(t-1)}[ \lambda (\lambda -1) (\lambda -t)\mu ^2 -\{ \theta _0  (\lambda -1) (\lambda -t)+\theta _1  \lambda (\lambda -t) \nonumber \\
\phantom{H=}{} +(\theta _t -1) \lambda (\lambda -1)\} \mu +\kappa _1 (\kappa _2 +1) (\lambda -t)], \label{eq:linP6}
\end{gather}
which we denote by $D_{y_1}(\theta _0, \theta _1, \theta _t, \theta _{\infty}; \lambda ,\mu )$.
This equation has regular singularities at $z=0,1,t,\lambda ,\infty$.
Exponents of the singularity $z=\lambda $ are $0$, $2$, and it is apparent (non-logarithmic) singularity.
Note that the dif\/ferential equations
\begin{gather}
\frac{d\lambda }{dt} =\frac{\partial H}{\partial \mu}, \qquad \frac{d\mu }{dt} =-\frac{\partial H}{\partial \lambda}
\label{eq:P6sys}
\end{gather}
describe the condition for monodromy preserving deformation of equation~(\ref{eq:y1DE}) with respect to the variable $t$.
By eliminating the variable $\mu$ in equation~(\ref{eq:P6sys}), we have the sixth Painlev\'e equation on the variable $\lambda $ (see equation~(\ref{eq:P6eqn})).
See~\cite{Oka1} on equations~(\ref{eq:y1DE}), (\ref{eq:linP6}) and (\ref{eq:P6sys}).

We consider realization of the Fuchsian system (equation~(\ref{eq:dYdzAzY1})) for the case $\lambda = 0,1,t,\infty$ in the appendix.

\section{The space of initial conditions for the sixth Painlev\'e equation\\ and Heun's equation} \label{sec:HeunSIC}

In this section, we introduce the space of initial conditions for the sixth Painlev\'e equation, restrict the variables of the space of initial conditions $ E(t)$ to certain lines, and we obtain Heun's equation.

The space of initial conditions was introduced by Okamoto \cite{Oka}, which is a suitable def\/ining variety for the set of solutions to the Painlev\'e system.
In \cite{ST}, Shioda and Takano studied the space of initial conditions further for the sixth Painlev\'e system (equation~(\ref{eq:P6sys})) to study roles of holomorphy on the Hamiltonian.
It was also constructed as a moduli space of parabolic connections by Inaba, Iwasaki and Saito \cite{IIS1,IIS2}.
Here we adopt the coordinate of initial coordinate by Shioda and Takano \cite{ST} (see also \cite{Wat}).
The space of initial condition $E(t)$ is def\/ined by patching six copies
\begin{gather}
  U_0=\{ (q_0, p_0) \}, \qquad U_1=\{ (q_1, p_1) \}, \qquad U_2=\{ (q_2, p_2) \} , \nonumber \\
  U_3=\{ (q_3, p_3) \}, \qquad U_4=\{ (q_4, p_4) \}, \qquad U_{\infty }=\{ (q_{\infty }, p_{\infty }) \}, \label{def:U01234i}
\end{gather}
of $\Cplx ^2$ for f\/ixed $(t; \theta_0, \theta _1,\theta _t,\theta _{\infty})$,
and the rule of patching is def\/ined by
\begin{alignat}{4}
& q_0 q_{\infty } =1, \qquad &&  q_0p_0+q_{\infty } p_{\infty } =-\kappa _1, \qquad  && (U_0 \cap U_{\infty}) , & \nonumber\\
& q_0p_0+q_1 p_1 =\theta _0 , \qquad && p_0p_1=1,   \qquad && (U_0 \cap U_1) ,& \nonumber \\
& (q_0-1)p_0+q_2 p_2 =\theta _1 , \qquad&& p_0p_2=1 ,  \qquad &&(U_0 \cap U_2) ,& \nonumber \\
& (q_0-t) p_0+q_3 p_3 =\theta _t , \qquad && p_0p_3=1 ,  \qquad & &(U_0 \cap U_3) ,& \nonumber \\
& q_{\infty } p_{\infty }+q_4 p_4 =1-\theta _{\infty } , \qquad&&  p_{\infty }p_4=1 ,   \qquad && (U_{\infty } \cap U_4) . & \label{eq:defrelSPI}
\end{alignat}
The variables $(\lambda ,\mu )$ of the sixth Painlev\'e system (see equation~(\ref{eq:P6sys})) are realized as $q_0=\lambda $, $p_0=\mu$ in $U_0$.

We def\/ine complex lines in the space of initial conditions as follows:
\begin{alignat}{3}
 & L_0 = \{ (0,p_0)\} \subset U_0, \qquad && L_1 = \{ (1,p_0)\} \subset U_0, & \nonumber\\
&  L_t = \{ (t,p_0)\} \subset U_0, \qquad && L_{\infty } = \{ (0,p_{\infty })\} \subset U_{\infty }, & \nonumber \\
 & L^*_0 = \{ (q_1,0)\} \subset U_1, \qquad && L^*_1 = \{ (q_2,0)\} \subset U_2,&  \nonumber \\
&  L^*_t = \{ (q_3,0)\} \subset U_3, \qquad && L^*_{\infty } = \{ (q_4,0)\} \subset U_{4}. & \label{eq:linesSIC}
\end{alignat}
Set
\begin{gather*}
  U_0^{q_0 \neq 0,1,t} = U_0 \setminus (L_0 \cup L_1 \cup L_t) .
\end{gather*}
Then the space of initial conditions $E(t)$ is a direct sum of the sets $U_0^{q_0 \neq 0,1,t}$, $L_0$, $L_1$, $L_t$, $L_{\infty}$, $L^*_0$, $L^*_1$, $L^*_t$, $L^*_{\infty}$.
If $(\lambda , \mu ) \in U_0^{q_0 \neq 0,1,t}$, then $\lambda \neq 0,1,t,\infty $ and equation~(\ref{eq:linP6}) has f\/ive regular singularities $\{ 0,1,t,\lambda , \infty \}$.

Although equation~(\ref{eq:wuA01t}) was considered on the set $ U_0^{q_0 \neq 0,1,t} $, we may consider realization of a~second-order dif\/ferential equation as equation~(\ref{eq:linP6}) on the space of initial conditions $E(t)$.
On the lines $L_0$, $L_1$, $L_t$, equation~(\ref{eq:linP6}) is realized by setting $\lambda =0 ,1,t$,
and the equation is written in the form of Heun's equation
\begin{gather}
  \frac{d^2y_1}{dz^2} + \left( \frac{-\theta _0}{z}+\frac{1-\theta _1}{z-1}+\frac{1-\theta _t}{z-t} \right)  \frac{dy_1}{dz}+  \frac{\kappa _1 (\kappa _2 +1) z+t\theta _0 \mu}{z(z-1)(z-t)} y_1=0 , \label{eq:HeunP60} \\
  \frac{d^2y_1}{dz^2} + \left( \frac{1-\theta _0}{z}+\frac{-\theta _1}{z-1}+\frac{1-\theta _t}{z-t} \right)  \frac{dy_1}{dz}+  \frac{\kappa _1 (\kappa _2 +1) (z-1)+(1-t)\theta _1 \mu }{z(z-1)(z-t)} y_1=0, \label{eq:HeunP61} \\
  \frac{d^2y_1}{dz^2} + \left( \frac{1-\theta _0}{z}+\frac{1-\theta _1}{z-1}+\frac{-\theta _t}{z-t} \right)  \frac{dy_1}{dz}+  \frac{\kappa _1 (\kappa _2 +1) (z-t)+t(t-1)\theta _t \mu}{z(z-1)(z-t)} y_1=0,  \label{eq:HeunP6t}
\end{gather}
respectively.
Note that if $\theta _0 \theta _1 \theta _t \neq 0$ then we can realize all values of accessory parameter as varying $\mu $.
For the case $\theta _0 \theta _1 \theta _t = 0$, we should consider other realizations.

To realize equation~(\ref{eq:linP6}) on the line $L^*_0$, we change the variables $(\lambda , \mu)$ into the ones $(q_1, p_1)$ on equation~(\ref{eq:linP6}) by applying relations $\lambda \mu +q_1 p_1 =\theta _0 , \; \mu p_1=1 $.
Then we have
\begin{gather*}
  \frac{d^2y_1}{dz^2} + \left( \frac{1-\theta _0}{z}+\frac{1-\theta _1}{z-1}+\frac{1-\theta _t}{z-t}-\frac{1}{z+p_1( p_1 q_1-\theta _0)} \right) \frac{dy_1}{dz} \\ %\label{eq:linP6q1p1} \\
 \phantom{\frac{d^2y_1}{dz^2}}{} + \frac{\kappa _1(\kappa _2 +1)z^2+ (t q_1-\theta _0 ( \theta _t+t \theta _1+p_1 \,{\rm pol}_1 )z+ ( p_1 q_1-\theta _0)(-t-p_1 \,{\rm pol}_2) }{z(z-1)(z-t)(z+p_1( p_1 q_1-\theta _0))}y_1 =0 ,\nonumber
\end{gather*}
where ${\rm pol}_1 $ and ${\rm pol}_2$ are polynomials in $p_1$, $q_1$, $t$, $\theta _0$, $\theta _1$, $\theta _t$, $\theta _{\infty }$.
By setting $p_1=0$, we obtain
\begin{gather}
  \frac{d^2y_1}{dz^2} +  \left( \frac{-\theta _0}{z}+\frac{1-\theta _1}{z-1}+\frac{1-\theta _t}{z-t} \right) \frac{dy_1}{dz}\nonumber\\
\phantom{\frac{d^2y_1}{dz^2}}{}
 + \frac{\kappa _1(\kappa _2 +1)z^2+ (t q_1-\theta _0 ( \theta _t+t \theta _1))z+ t \theta _0}{z^2(z-1)(z-t)}y_1 =0. \label{eq:linP6q1}
\end{gather}
Since the exponents of equation~(\ref{eq:linP6q1}) at $z=0$ are $1$ and $\theta _0$, we consider gauge-transformation $v_1 =z^{-1}y_1$ to obtain Heun's equation, and we have
\begin{gather}
  \frac{d^2v_1}{dz^2} + \left( \frac{2-\theta _0}{z}+\frac{1-\theta _1}{z-1}+\frac{1-\theta _t}{z-t} \right)  \frac{dv_1}{dz}+  \frac{(\kappa _1 +1)(\kappa _2 +2) z- q}{z(z-1)(z-t)} v_1=0,  \nonumber\\
  q=- tq_1 +(\theta _0-1)\{ t(\theta _1 -1) +\theta _t-1 \} . \label{eq:HeunP60v2}
\end{gather}
To realize the second-order Fuchsian equation on the line $L^*_1$, we change the variables $(\lambda , \mu)$ into the ones $(q_2, p_2)$, substitute $p_2=0 $ into equation~(\ref{eq:linP6}) and set $v_1 =(z-1)^{-1}y_1$.
Then $v_1 $ satisf\/ies the following equation;
\begin{gather}
  \frac{d^2v_1}{dz^2} + \left( \frac{1-\theta _0}{z}+\frac{2-\theta _1}{z-1}+\frac{1-\theta _t}{z-t} \right)  \frac{dv_1}{dz}+  \frac{(\kappa _1 +1)(\kappa _2 +2) (z-1)- q}{z(z-1)(z-t)} v_1=0,  \nonumber\\
  q=(t-1)q_2 -(\theta _1 -1) \{ (1-t) (\theta _0-1)+\theta _t -1 \}. \label{eq:HeunP61v2}
\end{gather}
The second-order Fuchsian equation on the line $L^*_t$ is realized as
\begin{gather}
  \frac{d^2v_1}{dz^2} + \left( \frac{1-\theta _0}{z}+\frac{1-\theta _1}{z-1}+\frac{2-\theta _t}{z-t} \right)  \frac{dv_1}{dz}+  \frac{(\kappa _1 +1)(\kappa _2 +2) (z-t)- q}{z(z-1)(z-t)} v_1=0,  \nonumber\\
  q=t(1-t)q_3 -(\theta _t-1)((t-1)(\theta _0 -1)+t(\theta _1-1)), \label{eq:HeunP6tv2}
\end{gather}
by setting $p_3=0$ and $v_1 =(z-t)^{-1}y_1$.

We investigate equation~(\ref{eq:linP6}) on the line $L_{\infty }$.
We change the variables $(\lambda , \mu)$ into the ones $(q_{\infty }, p_{\infty })$ on equation~(\ref{eq:linP6}) by applying relations $\lambda q_{\infty } =1, \; \lambda \mu +q_{\infty } p_{\infty } =-\kappa _1$, and substitute $q_{\infty }=0 $.
Then we have
\begin{gather*}
  \frac{d^2y_1}{dz^2} + \left( \frac{1-\theta _0}{z}+\frac{1-\theta _1}{z-1}+\frac{1-\theta _t}{z-t} \right)  \frac{dy_1}{dz}+  \frac{\kappa _1 (\kappa _2 +2) z- q}{z(z-1)(z-t)} y_1=0, \nonumber \\
  q=(\theta _{\infty }-1) p_{\infty }  +\kappa _1 ( t(\kappa _2  +\theta _t +1)+\kappa _2 +\theta _1 +1).%\label{eq:HeunP6infv1}
\end{gather*}
Note that the exponents at $z=\infty $ are $\kappa _1$ and $\kappa _2+2$.

To realize equation~(\ref{eq:linP6}) on the line $L^*_{\infty }$, we change the variables $(\lambda , \mu)$ into the ones $(q_{4 }, p_{4})$ on equation~(\ref{eq:linP6}) by applying relations $\lambda q_{\infty } =1$, $\lambda \mu +q_{\infty } p_{\infty } =-\kappa _1$, $q_{\infty } p_{\infty }+q_4 p_4 =1-\theta _{\infty }$, $p_{\infty }p_4=1$, substitute $p_{4 }=0 $. We obtain
\begin{gather}
  \frac{d^2y_1}{dz^2} + \left( \frac{1-\theta _0}{z}+\frac{1-\theta _1}{z-1}+\frac{1-\theta _t}{z-t} \right)  \frac{dy_1}{dz}+  \frac{(\kappa _1 +1)(\kappa _2 +1) z- q}{z(z-1)(z-t)} y_1=0,  \nonumber\\
  q=-q_{4 } + (\kappa _2 +1)(t(\kappa _1+ \theta _t )+\kappa _1+\theta _1).\label{eq:HeunP6infv2}
\end{gather}
The exponents at $z=\infty $ are $\kappa _1+1 $ and $\kappa _2+1$, which are dif\/ferent from the case of the line $L_{\infty }$.

The Fuchsian system $D_Y(\theta _0, \theta _1, \theta _t, \theta _{\infty}; \lambda ,\mu ;k )$ is originally def\/ined on the set $ U_0^{q_0 \neq 0,1,t} $.
We try to consider realization of Fuchsian system (equation~(\ref{eq:dYdzAzY1})) on the lines $L_0$, $L^*_0$, $L_1$, $L^*_1$, $L_t$, $L^*_t$, $L_{\infty}$, $L^*_{\infty}$ in the appendix.

\section{Middle convolution} \label{sec:MC}

First, we review an algebraic analogue of Katz' middle convolution functor developed by Dettweiler and Reiter \cite{DR1,DR2}, which we restrict to the present setting.
Let $A_0$, $A_1$, $A_t$ be matrices in~$\Cplx ^{2\times 2}$.
For $\nu \in \Cplx$, we def\/ine the convolution matrices $B_0, B_1, B_t \in \Cplx ^{6\times 6}$ as follows:
\begin{gather}
  B_0=
\left(
\begin{array}{ccc}
A_0 +\nu & A_1 & A_t \\
0 & 0 & 0 \\
0 & 0 & 0
\end{array}
\right) , \qquad
B_1=
\left(
\begin{array}{ccc}
0 & 0 & 0 \\
A_0 & A_1 +\nu & A_t \\
0 & 0 & 0
\end{array}
\right) ,   \nonumber \\
  B_t=
\left(
\begin{array}{ccc}
0 & 0 & 0 \\
0 & 0 & 0 \\
A_0 & A_1 & A_t +\nu
\end{array}
\right) . \label{eq:Bdef}
\end{gather}
Let $z \in \Cplx \setminus \{0,1,t\}$, $\gamma _p$ $(p\in \Cplx )$ be a cycle in $\Cplx \setminus \{ 0,1,t,z \}$ turning the point $w=p$ anti-clockwise whose f\/ixed base point is $o \in \Cplx \setminus \{0,1,t ,z \}$,
and $[\gamma _p , \gamma _{p'}] = \gamma _p  \gamma _{p'} \gamma ^{-1}_p \gamma ^{-1}_{p'} $ be the Pochhammer contour.

\begin{proposition}[\cite{DR2}] \label{prop:DRintegrepr}
Assume that $Y= {}^t (y_1(z), y_2(z))$ is a solution to the system of differential equations
\begin{gather*}
\frac{dY}{dz}=\left( \frac{A_0}{z}+\frac{A_1}{z-1}+\frac{A_t}{z-t} \right) Y.
%\label{eq:dYdzAzY0}
\end{gather*}
For $p \in \{ 0,1,t,\infty \}$, the function
\begin{gather*}
U = \left(
\begin{array}{l}
\displaystyle \int _{[\gamma _{z} ,\gamma _p]}w^{-1}y_1(w) (z-w)^{\nu } dw \vspace{1mm}\\
\displaystyle \int _{[\gamma _{z} ,\gamma _p]}w^{-1}y_2(w) (z-w)^{\nu } dw \vspace{1mm}\\
\displaystyle \int _{[\gamma _{z} ,\gamma _p]}(w-1)^{-1}y_1(w) (z-w)^{\nu } dw \vspace{1mm}\\
\displaystyle \int _{[\gamma _{z} ,\gamma _p]}(w-1)^{-1}y_2(w) (z-w)^{\nu } dw \vspace{1mm}\\
\displaystyle \int _{[\gamma _{z} ,\gamma _p]}(w-t)^{-1}y_1(w) (z-w)^{\nu } dw \vspace{1mm}\\
\displaystyle \int _{[\gamma _{z} ,\gamma _p]}(w-t)^{-1}y_2(w) (z-w)^{\nu } dw
\end{array}
\right) ,% \label{eq:integrepU}
\end{gather*}
satisfies the system of differential equations
\begin{gather}
\frac{dU}{dz}=\left( \frac{B_0}{z}+\frac{B_1}{z-1}+\frac{B_t}{z-t} \right) U.
\label{eq:dYdzBzU}
\end{gather}
\end{proposition}
We set
\begin{gather}
  {\mathcal L}_0= \left(
\begin{array}{c}
\mbox{Ker}(A_0) \\
0 \\
0
\end{array}
\right) , \qquad
{\mathcal L}_1= \left(
\begin{array}{c}
0\\
\mbox{Ker}(A_1) \\
0
\end{array}
\right) , \qquad
{\mathcal L}_t= \left(
\begin{array}{c}
0 \\
0 \\
\mbox{Ker}(A_t)
\end{array}
\right) ,\nonumber  \\
  {\mathcal L} ={\mathcal L}_0 \oplus {\mathcal L}_1 \oplus {\mathcal L}_t, \qquad
{\mathcal K} = \mbox{Ker}(B_0) \cap \mbox{Ker}(B_1) \cap  \mbox{Ker}(B_t) , \label{eq:KL}
\end{gather}
where ${\mathcal L}_0$, ${\mathcal L}_1$, ${\mathcal L}_t$, ${\mathcal K} \subset \Cplx ^6 $ and $0$ in equation~(\ref{eq:KL}) means the zero vector in $\Cplx ^2$.
 We f\/ix an isomorphism between $\Cplx ^6 /({\mathcal K}+ {\mathcal L})$ and $\Cplx ^m$ for some $m$.
A tuple of matrices $mc_{\nu } (A) =(\tilde{B}_0, \tilde{B}_1, \tilde{B}_t)$, where $\tilde{B}_p$ $(p=0,1,t)$ is induced by the action of $B_p$ on $\Cplx ^m \simeq \Cplx ^6 /({\mathcal K}+ {\mathcal L})$, is called an additive version of the middle convolution of $(A_0,A_1,A_t)$ with the parameter $\nu$.
Let~$A_0$,~$A_1$,~$A_t$ be the matrices def\/ined by equation~(\ref{eq:A0A1AtP}).
Then it is shown that if $\nu =0, \kappa _1, \kappa _2$ (resp.\ $\nu \neq 0, \kappa _1, \kappa _2$) then $\dim \Cplx ^6 /({\mathcal K}+ {\mathcal L}) =2$ (resp.\ $\dim \Cplx ^6 /({\mathcal K}+ {\mathcal L}) =3$).
If $\nu =0$, then the middle convolution is identity (see~\cite{DR2}).
Hence the middle convolutions for two cases $\nu = \kappa _1,\kappa _2$ may lead to non-trivial transformations on the $2 \times 2$ Fuchsian system with four singularities $\{ 0,1,t,\infty \}$.
Filipuk \cite{Fil} obtained that the middle convolution for the case $\nu =\kappa _1$ induce an Okamoto's transformation of the sixth Painlev\'e system.

We now calculate explicitly the Fuchsian system of dif\/ferential equations determined by the middle convolution for the case $\nu =\kappa _2$.
Note that the following calculation is analogous to the one in \cite{TakI} for the case $\nu =\kappa _1$.
If $\nu = \kappa _2$, then the spaces ${\mathcal L}_0$, ${\mathcal L}_1$, ${\mathcal L}_t$, ${\mathcal K}$ are written as
\begin{gather*}
  {\mathcal L}_0= \Cplx \left(\!\!
\begin{array}{c}
w_0 \\
u_0+\theta _0 \\
0\\
0\\
0\\
0
\end{array}\!\!
\right) , \quad \ \
{\mathcal L}_1= \Cplx \left(\!\!
\begin{array}{c}
0\\
0\\
w_1 \\
u_1+\theta _1 \\
0\\
0
\end{array}\!\!
\right) , \quad \ \
{\mathcal L}_t= \Cplx \left(\!\!
\begin{array}{c}
0\\
0\\
0\\
0\\
w_t \\
u_t+\theta _t
\end{array}\!\!
\right) ,\quad \ \
{\mathcal K}= \Cplx \left(\!\!
\begin{array}{c}
0\\
1\\
0\\
1\\
0\\
1
\end{array}\!\!
\right) .
\end{gather*}
Set
\begin{gather}
  S= \left(
\begin{array}{cccccc}
0 & 0 & 0 & w_0 & 0 & 0\\
0 & 0 & 1 &u_0+\theta _0& 0 & 0 \\
0 & 0 & 0 &0 & w_1 & 0 \\
s_{41} & s_{42} & 1 & 0 & u_1+\theta _1& 0 \\
0 & 0 & 0 &0& 0 &  w_t \\
s_{61} & s_{62} & 1 &0 & 0 & u_t+\theta _t
\end{array}
\right) ,  \nonumber\\
s_{41}= \frac{\mu (\lambda -t)+\kappa _1}{k \kappa _1 }, \qquad
s_{61}= \frac{t(\mu (\lambda -1)+\kappa _1)}{k \kappa _1 }, \qquad
s_{42}= \frac{\tilde{\lambda } -\lambda}{\lambda (\lambda -1) \kappa _2} ,\nonumber\\
s_{62}= \frac{t(\tilde{\lambda } -\lambda )}{\lambda (\lambda -t) \kappa _2},\qquad
 \tilde{\lambda }=\lambda -\frac{\kappa _2}{\mu -\frac{\theta _0}{\lambda }-\frac{\theta _1}{(\lambda -1)}-\frac{\theta _t}{(\lambda -t)}}, \label{eq:matS}
\end{gather}
and $\tilde{U}= S^{-1} U$,  where $U$ is a solution to equation~(\ref{eq:dYdzBzU}).
Then $\det U =k ^2(\tilde{\lambda }- \lambda )/(t(1-t)\kappa _2)$ and $\tilde{U}$ satisf\/ies
\begin{gather*}
\frac{d\tilde{U}}{dz} = \left(
\begin{array}{cccccc}
b_{11}(z) & b_{12}(z) & 0 & 0 & 0 & 0\\
b_{21}(z) & b_{22}(z) & 0 & 0 & 0 & 0\vspace{1mm}\\
\frac{-(u_0+\theta _0)\theta _{\infty }t}{k\kappa _1 \lambda z} & \frac{\tilde{\lambda }}{\lambda z} & 0 & 0 & 0 & 0 \vspace{1mm}\\
\frac{t}{k\lambda z} & 0 & 0 & \frac{\kappa _2}{z} & 0 & 0\vspace{1mm}\\
\frac{1-t}{k(\lambda -1)(z-1)} & 0 & 0 & 0 & \frac{\kappa _2}{z-1} & 0\vspace{1mm}\\
\frac{t(1-t)}{k(\lambda -t)(z-t)} & 0 & 0 & 0 & 0 & \frac{\kappa _2}{z-t}
\end{array}
\right) \tilde{U},
\end{gather*}
where $b_{11}(z), \dots , b_{22}(z)$ are calculated such that the system of dif\/ferential equation
\begin{gather*}
  \frac{d\tilde{Y}}{dz}=
\left(
\begin{array}{ll}
b_{11}(z) & b_{12}(z) \\
b_{21}(z) & b_{22}(z)
\end{array}
\right) \tilde{Y}, \qquad
\tilde{Y}=
\left(
\begin{array}{ll}
\tilde{u}_1(z) \\
\tilde{u}_2(z)
\end{array}
\right)  ,
%\label{eq:dtYdzBztY0}
\end{gather*}
coincides with the Fuchsian system $D_Y(\tilde{\theta }_0, \tilde{\theta }_1, \tilde{\theta }_t, \tilde{\theta }_{\infty}; \tilde{\lambda },\tilde{\mu };\tilde{k} )$ (see equation~(\ref{eq:dYdzAzY})), where
\begin{gather}
 \tilde{\theta }_0 = \frac{\theta _0 -\theta _1 -\theta _t-\theta _{\infty}}{2}, \qquad \tilde{\theta }_1 = \frac{-\theta _0 +\theta _1 -\theta _t-\theta _{\infty}}{2}, \qquad \tilde{\theta }_t = \frac{-\theta _0 -\theta _1 +\theta _t-\theta _{\infty}}{2}, \nonumber \\
   \tilde{\theta }_{\infty} = \frac{-\theta _0-\theta _1 -\theta _t+\theta _{\infty}}{2}, \qquad \tilde{\lambda} =\lambda -\frac{\kappa _2}{\mu -\frac{\theta _0}{\lambda }-\frac{\theta _1}{\lambda -1}-\frac{\theta _t}{\lambda -t}} ,\nonumber \\
  \tilde{\mu} =\frac{\kappa _2 +\theta _0}{\tilde{\lambda }} +\frac{\kappa _2 +\theta _1}{\tilde{\lambda }-1} +\frac{\kappa _2 +\theta _t}{\tilde{\lambda }-t} +\frac{\kappa _2 }{\lambda - \tilde{\lambda }} ,\qquad \tilde{k}= k . \label{eq:tt0tt1ttt}
\end{gather}
The functions $\tilde{u}_1(z)$ and $\tilde{u}_2(z)$ are expressed as
\begin{gather}
  \tilde{u}_1(z) = (u_0+\theta _0) u_1(z) - \frac{k\lambda }{t} u_2 (z) + (u_1+\theta _1) u_3(z)  \nonumber\\
\phantom{\tilde{u}_1(z) =}{}
 + \frac{k(\lambda -1)}{t-1 } u_4 (z)+ (u_t+\theta _t) u_5(z) + \frac{k(\lambda -t)  }{t(1-t)} u_6 (z), \nonumber \\
  \tilde{u}_2(z) = \frac{\kappa _2 \lambda (\lambda -1) (\lambda -t)}{\kappa _1(\lambda - \tilde{\lambda })} \left( \frac{(\lambda \mu +\kappa _1 ) (u_0 + \theta _0 )}{k \lambda } u_1(z) - \frac{\lambda \mu +\kappa _1 }{t} u_2 (z) \right. \nonumber \\
\phantom{\tilde{u}_2(z) =}{}   + \frac{((\lambda -1)\mu +\kappa _1 ) (u_1 + \theta _1 )}{k (\lambda -1)} u_3(z)+ \frac{(\lambda -1)\mu +\kappa _1 }{t-1} u_4 (z) \nonumber \\
  \left.\phantom{\tilde{u}_2(z) =}{}
 +\frac{((\lambda -t)\mu +\kappa _1 ) (u_t + \theta _t )}{k (\lambda -t)} u_5(z) +\frac{(\lambda -t)\mu +\kappa _1 }{t (1-t)} u_6 (z) \right) .\label{eq:y1tu16}
\end{gather}
Combining Proposition \ref{prop:DRintegrepr} with equation~(\ref{eq:y1tu16}) and setting $\tilde{y}_1(z)=\tilde{u}_1(z)$, $\tilde{y}_2(z)=\tilde{u}_2(z)$, we have the following theorem by means of a straightforward calculation:

\begin{theorem} \label{thm:zinterep}
Set $\kappa _1= (\theta _{\infty } -\theta _0 -\theta _1 -\theta _t)/2$ and $\kappa _2= -(\theta _{\infty } +\theta _0 +\theta _1 +\theta _t)/2$.
If $y_1(z)$ is a~solution to the Fuchsian equation $D_{y_1}(\theta _0, \theta _1, \theta _t,\theta _{\infty}; \lambda ,\mu )$, then the function $\tilde{Y} = {}^t (\tilde{y}_1 (z) , \tilde{y}_2 (z) )$ defined by
\begin{gather}
  \tilde{y}_1 (z) = \int _{[\gamma _{z} ,\gamma _p]}  \frac{dy_1(w)}{dw}  (z-w)^{\kappa _2 } dw, \label{eq:yt1zintrep} \\
  \tilde{y}_2 (z)= \frac{\kappa _2 \lambda (\lambda -1) (\lambda -t)}{k(\lambda -\tilde{\lambda })} \!\int _{[\gamma _{z} ,\gamma _p]}\! \left\{ \left( \frac{dy_1(w)}{dw} -\mu y_1(w)  \right)\! \frac{1}{\lambda -w} +\frac{\mu}{\kappa _1}\frac{dy_1(w)}{dw} \right\} (z-w)^{\kappa _2 } dw,  \nonumber
\end{gather}
satisfies the Fuchsian system $D_Y(\kappa _2 + \theta _0, \kappa _2 +\theta _1, \kappa _2 +\theta _t, \kappa _2 +\theta _{\infty} ; \tilde{\lambda } ,\tilde{\mu };k )$ for $p \in \{ 0,1,t,\infty \}$, where
\begin{gather}
  \tilde{\lambda} =\lambda -\frac{\kappa _2}{\mu -\frac{\theta _0}{\lambda }-\frac{\theta _1}{\lambda -1}-\frac{\theta _t}{\lambda -t}} ,\qquad  \tilde{\mu} =\frac{\kappa _2 +\theta _0}{\tilde{\lambda }} +\frac{\kappa _2 +\theta _1}{\tilde{\lambda }-1} +\frac{\kappa _2 +\theta _t}{\tilde{\lambda }-t} +\frac{\kappa _2 }{\lambda - \tilde{\lambda }} . \label{eq:tlamtmu}
\end{gather}
\end{theorem}
Since
\begin{gather}
  0= \int _{[\gamma _{z} ,\gamma _p]}  \frac{d}{dw} \left(  y_1(w)(z-w)^{\kappa _2 } \right) dw \nonumber \\
\phantom{0}{} = \int _{[\gamma _{z} ,\gamma _p]}  \frac{dy_1(w)}{dw}  (z-w)^{\kappa _2 } dw + \kappa _2 \int _{[\gamma _{z} ,\gamma _p]} y_1(w) (z-w)^{\kappa _2 -1} dw, \label{eq:intbyparts0}
\end{gather}
we have
\begin{proposition}[\cite{Nov}] \label{prop:Nov}
If $y_1(z)$ is a solution to $D_{y_1}(\theta _0, \theta _1, \theta _t,\theta _{\infty}; \lambda ,\mu )$,
then the function
\begin{gather}
\tilde{y} (z)  = \int _{[\gamma _{z} ,\gamma _p]} y_1(w) (z-w)^{\kappa _2 -1} dw,
\label{eq:inttransDy1}
\end{gather}
satisfies $D_{y_1}(\kappa _2 + \theta _0, \kappa _2 +\theta _1, \kappa _2 +\theta _t, \kappa _2 +\theta _{\infty} ; \tilde{\lambda } ,\tilde{\mu } )$ for $p \in \{ 0,1,t,\infty \}$, where $\tilde{\lambda }$ and $\tilde{\mu } $ are defined in equation~\eqref{eq:tlamtmu}.
\end{proposition}
Note that this proposition was obtained by Novikov \cite{Nov} by another method.
Kazakov and Slavyanov~\cite{KS1} essentially obtained this proposition by investigating Euler transformation of $2\times 2$ Fuchsian systems with singularities $\{0,1,t,\infty \}$ which are realized dif\/ferently from $D_{Y}(\theta _0, \theta _1, \theta _t,\theta _{\infty}; \lambda ,\mu ,k)$.

Let us recall the symmetry of the sixth Painlev\'e equation.
It was essentially established by Okamoto~\cite{Oka0} that the sixth Painlev\'e equation has symmetry of the af\/f\/ine Weyl group $W(D_4^{(1)})$.
More precisely, the sixth Painlev\'e system is invariant under the following transformations, which are involutive and satisfy Coxeter relations attached to the Dynkin diagram of type $D_4^{(1)}$, i.e.\ $(s_i )^2=1$ $(i=0,1,2,3,4)$, $s_j s_k=s_k s_j $ $(j,k \in \{ 0,1,3,4 \})$, $s_j s_2 s_j= s_2 s_j s_2$ $(j=0,1,3,4)$:

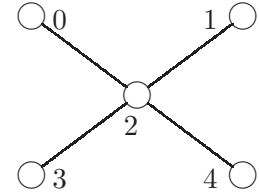
\begin{tabular}{|c||cccc|cc|c|}\hline
{} & $\theta _t $&$\theta _{\infty} $&$\theta _1 $&$\theta _0 $&$ \lambda  $&$ \mu  $&$ t$\\
\hline
$s_0$&$ -\theta _t $&$\theta _{\infty} $&$\theta _1 $&$\theta _0 $&$ \lambda  $&$ \mu -\frac{\theta _t}{\lambda -t} $&$ t$\\
$s_1$&$\theta _t $&$ 2-\theta _{\infty} $&$\theta _1 $&$\theta _0 $&$ \lambda  $&$ \mu $&$ t$\\
$s_2$&$\kappa _1+ \theta _t $&$-\kappa _2 $&$\kappa _1+\theta _1 $&$\kappa _1+\theta _0 $&$ \lambda +\frac{\kappa _1}{\mu } $&$ \mu  $&$ t$\\
$s_3$&$\theta _t $&$\theta _{\infty} $&$ -\theta _1 $&$\theta _0 $&$ \lambda  $&$ \mu -\frac{\theta _1}{\lambda -1} $&$ t$\\
$s_4 $&$\theta _t $&$\theta _{\infty} $&$\theta _1 $&$ -\theta _0 $&$ \lambda  $&$ \mu -\frac{\theta _0}{\lambda } $&$ t$\\
\hline
\end{tabular}
\begin{picture}(100,55)(50,45)
\put(70,80){\circle{10}}
\put(70,20){\circle{10}}
\put(150,80){\circle{10}}
\put(150,20){\circle{10}}
\put(110,50){\circle{10}}
\qbezier(74,77)(90,65)(106,53)
\qbezier(74,23)(90,35)(106,47)
\qbezier(146,77)(130,65)(114,53)
\qbezier(146,23)(130,35)(114,47)
\put(78,75){$0$}
\put(78,15){$3$}
\put(135,75){$1$}
\put(135,15){$4$}
\put(105,35){$2$}
\end{picture}

\vspace{1mm}

The map $(\theta _0, \theta _1, \theta _t,\theta _{\infty}; \lambda ,\mu ) \mapsto (\tilde{\theta }_0, \tilde{\theta }_1, \tilde{\theta }_t, \tilde{\theta }_{\infty}; \tilde{\lambda },\tilde{\mu })$ determined by equation~(\ref{eq:tt0tt1ttt}) coincides with the composition map $s_0s_3s_4s_2s_0s_3s_4$, because
\begin{gather*}
  (\theta _0, \theta _1, \theta _t,\theta _{\infty}; \lambda ,\mu ) \mathop{\mapsto }^{s_0s_3s_4}
\textstyle \big(-\theta _0, -\theta _1, -\theta _t,\theta _{\infty}; \lambda ,\mu -\frac{\theta _0}{\lambda }-\frac{\theta _1}{\lambda -1}-\frac{\theta _t}{\lambda -t}\big) \\
\phantom{(\theta _0, \theta _1, \theta _t,\theta _{\infty}; \lambda ,\mu )}{}
  \mathop{\mapsto }^{s_2} \textstyle \big(-\kappa _2 -\theta _0, -\kappa _2 -\theta _1, -\kappa _2 -\theta _t, \kappa _1 ; \tilde{\lambda },\mu -\frac{\theta _0}{\lambda }-\frac{\theta _1}{\lambda -1}-\frac{\theta _t}{\lambda -t}\big) \nonumber \\
\phantom{(\theta _0, \theta _1, \theta _t,\theta _{\infty}; \lambda ,\mu )}{}
  \mathop{\mapsto }^{s_0s_3s_4} \textstyle \big(\kappa _2 +\theta _0, \kappa _2 +\theta _1, \kappa _2 +\theta _t, \kappa _2 + \theta _{\infty} ; \tilde{\lambda },\tilde{\mu}\big) . \nonumber
\end{gather*}
Therefore, if we know a solution to the Fuchsian system $D_Y(\theta _0, \theta _1, \theta _t,\theta _{\infty}; \lambda ,\mu ;k )$, then we have integral representations of solutions to the Fuchsian system $D_Y(\tilde{\theta }_0, \tilde{\theta }_1, \tilde{\theta }_t, \tilde{\theta }_{\infty}; \tilde{\lambda },\tilde{\mu };k )$ obtained by the transformation $s_0s_3s_4s_2s_4s_3s_0$.
Note that the transformations $s_i$ $(i=0,1,2,3,4)$ are extended to isomorphisms of the space of initial conditions $E(t)$.

We recall the middle convolution for the case $\nu =\kappa _1$.
\begin{proposition}[\protect{\cite[Proposition 3.2]{TakI}}] \label{prop:zinterepk1}
If $Y={}^t ( y_1(z), y_2(z))$ is a solution to the Fuchsian system $D_Y(\theta _0, \theta _1, \theta _t,\theta _{\infty}; \lambda ,\mu ;k )$ $($see equation~\eqref{eq:dYdzAzY}$)$, then the function $\tilde{Y} = {}^t (\tilde{y}_1 (z) , \tilde{y}_2 (z) )$ defined by
\begin{gather}
  \tilde{y}_1 (z)= \int _{[\gamma _{z} ,\gamma _p]} \left\{ \kappa _1 y_1(w) + (w- \tilde{\lambda})\frac{dy_1(w)}{dw}  \right\} \frac{(z-w)^{\kappa _1 }}{w-\lambda } dw, \nonumber \\
  \tilde{y}_2 (z) = \frac{-\theta _{\infty}}{\kappa _2 }\int _{[\gamma _{z} ,\gamma _p]}  \frac{dy_2(w)}{dw}  (z-w)^{\kappa _1 } dw,  \label{eq:yt1zintrepk1}
\end{gather}
satisfies the Fuchsian system $D_Y(\kappa _1 + \theta _0, \kappa _1 +\theta _1, \kappa _1 +\theta _t, -\kappa _2 ; \lambda +\kappa _1 /\mu ,\mu ;k )$ for $p \in \{ 0,1,t,\infty \}$.
\end{proposition}
The parameters $(\kappa _1 + \theta _0, \kappa _1 +\theta _1, \kappa _1 +\theta _t, -\kappa _2 ; \lambda +\kappa _1 /\mu ,\mu )$ are obtained from the parameters $(\theta _0, \theta _1, \theta _t,\theta _{\infty}; \lambda ,\mu )$ by applying the transformation $s_2$.
Note that the relationship the transformation $s_2$ was obtained by Filipuk~\cite{Fil} explicitly (see also~\cite{HF}), and Boalch~\cite{Boa} and Dettweiler and Reiter~\cite{DR3} also obtained results on the symmetry of the sixth Painlev\'e equation and the middle convolution.

\section{Middle convolution, integral transformations\\ of Heun's equation and the space of initial conditions} \label{sec:MCHeun}

In this section, we investigating relationship among the middle convolution, integral transformations of Heun's equation and the space of initial conditions.

Kazakov and Slavyanov established an integral transformation on solutions to Heun's equation in \cite{KS}, which we express in a slightly dif\/ferent form.
\begin{theorem}[\cite{KS}] \label{thm:Heunint}
Set
\begin{gather}
 (\eta -\alpha )(\eta -\beta )=0 , \qquad \gamma '=\gamma +1 -\eta, \qquad \delta' =\delta +1-\eta , \qquad \epsilon '=\epsilon +1-\eta , \nonumber \\
 \{ \alpha ' , \beta ' \} = \{ 2-\eta , -2\eta + \alpha +\beta +1 \} , \nonumber \\
 q'=q+(1-\eta )(\epsilon +\delta t+(\gamma -\eta ) (t+1)). \label{eq:mualbe}
\end{gather}
Let $v(w)$ be a solution to
\begin{gather}
\frac{d^2v}{dw^2} + \left( \frac{\gamma '}{w}+\frac{\delta '}{w-1}+\frac{\epsilon '}{w-t}\right) \frac{dv}{dw} +\frac{\alpha ' \beta ' w -q'}{w(w-1)(w-t)} v=0.
\label{Heun01}
\end{gather}
Then the function
\begin{gather*}
  y(z)=\int _{[\gamma _z ,\gamma _p]} v(w) (z-w)^{-\eta } dw
\end{gather*}
is a solution to
\begin{gather}
\frac{d^2y}{dz^2} + \left( \frac{\gamma}{z}+\frac{\delta }{z-1}+\frac{\epsilon}{z-t}\right) \frac{dy}{dz} +\frac{\alpha \beta z -q}{z(z-1)(z-t)} y=0,
\label{Heun02}
\end{gather}
for $p \in \{ 0,1,t,\infty \}$.
\end{theorem}

Here we derive Theorem~\ref{thm:Heunint} by considering the limit $\lambda \rightarrow 0$ in Proposition~\ref{prop:Nov}.
Let us recall notations in Proposition~\ref{prop:Nov}.
We consider the limit $\lambda \rightarrow 0$ while f\/ixing $\mu$ for the case $\theta _0 \neq 0$ and $\theta _0 +\kappa _2 \neq 0$.
Then we have $\tilde{\lambda } \rightarrow 0$ and $\tilde{\mu } \rightarrow (t \theta _0 \mu+\kappa _2 (t(\kappa _1 +\theta _t)+\kappa _1+\theta _1))/(t(\kappa _2 +\theta _0))$.
Hence it follows from Proposition \ref{prop:Nov} and equation~(\ref{eq:HeunP60}) that, if $y(z)$ satisf\/ies
\begin{gather}
  \frac{d^2y(z)}{dz^2} + \left( \frac{-\theta _0}{z}+\frac{1-\theta _1}{z-1}+\frac{1-\theta _t}{z-t} \right)  \frac{dy(z)}{dz}  + \frac{\kappa _1(\kappa _2 +1)z+t\theta _0 \mu }{z(z-1)(z-t)} y(z)=0,  \label{eq:linP6H1}
\end{gather}
then the function
\begin{gather}
\tilde{y} (z)  = \int _{[\gamma _{z} ,\gamma _p]} y(w) (z-w)^{\kappa _2 -1} dw, \label{eq:linP6H12inttrans}
\end{gather}
satisf\/ies
\begin{gather}
  \frac{d^2\tilde{y}(z)}{dz^2} + \left( \frac{-\kappa _2-\theta _0 }{z}+\frac{1-\kappa _2-\theta _1}{z-1}+\frac{1-\kappa _2-\theta _t}{z-t} \right)  \frac{d\tilde{y}(z)}{dz}  \nonumber  \\
\phantom{\frac{d^2\tilde{y}(z)}{dz^2}}{}  + \left( \frac{\theta _{\infty }(1- \kappa _2 )z+ t(\kappa _2 + \theta _0 )\frac{t \theta _0 \mu+\kappa _2 (t(\kappa _1 +\theta _t)+\kappa _1+\theta _1))}{t(\kappa _2 +\theta _0)} }{z(z-1)(z-t)}  \right) \tilde{y}(z)=0.   \label{eq:linP6H2}
\end{gather}
By setting $\gamma =-\kappa _2 -\theta _0 $, $\delta = 1-\kappa _2-\theta _1$, $\epsilon = 1-\kappa _2-\theta _t$, $\alpha = \eta =1-\kappa _2$, $\beta =\theta _{\infty }$, $q= -\{ t \theta _0 \mu+\kappa _2 (t(\kappa _1 +\theta _t)+\kappa _1+\theta _1)) \}$ and comparing with the standard form of Heun's equation (equation~(\ref{eq:Heun})), we recover Theorem~\ref{thm:Heunint}.
Note that we can obtain the formula corresponding to the case $\theta _0=0$ (resp.\ $\theta _0 +\kappa _2 = 0$) by considering the limit $\theta _0 \rightarrow 0$ (resp.\ $\theta _0 +\kappa _2 \rightarrow 0$).

The limit $\lambda \rightarrow 0$ while f\/ixing $\mu$ implies the restriction of the coordinate $(\lambda ,\mu )$ to the line $L_0$ in the space of initial conditions $E(t)$, and the line $L_0$ with the parameter $(\theta _0, \theta _1, \theta _t,\theta _{\infty}; \lambda ,\mu )$ is mapped to the line $L_0$ in the space of initial conditions with the parameter $(\kappa _2 + \theta _0, \kappa _2 +\theta _1, \kappa _2 +\theta _t$, $\kappa _2 +\theta _{\infty} ; \tilde{\lambda } ,\tilde{\mu } )$ where $\tilde{\lambda }$ and $\tilde{\mu } $ are def\/ined in equation~(\ref{eq:tlamtmu}), because $\tilde{\lambda } \rightarrow 0$ and $\tilde{\mu }$ converges by the limit.
It follows from equations~(\ref{eq:linP6H1}), (\ref{eq:linP6H12inttrans}), (\ref{eq:linP6H2}) that the integral transformation in Proposition \ref{prop:Nov} reproduces the integral transformation on Heun's equations in Theorem \ref{thm:Heunint} by restricting to the line $L_0$.
We can also establish that the line $L_1$ (resp. $L_t$) in the space of initial conditions with the parameter $(\theta _0, \theta _1, \theta _t,\theta _{\infty}; \lambda ,\mu )$ is mapped to the line $L_1$ (resp. $L_t$) with the parameter $(\kappa _2 + \theta _0, \kappa _2 +\theta _1, \kappa _2 +\theta _t, \kappa _2 +\theta _{\infty} ; \tilde{\lambda } ,\tilde{\mu } )$ by taking the limit $\lambda \rightarrow 1$ (resp. $\lambda \rightarrow t$),
and the integral transformation in Proposition \ref{prop:Nov} reproduces the integral transformation on Heun's equations in Theorem \ref{thm:Heunint}.
We discuss the restriction of the map $(\theta _0, \theta _1, \theta _t,\theta _{\infty}; \lambda ,\mu ) \mapsto (\kappa _2 + \theta _0, \kappa _2 +\theta _1, \kappa _2 +\theta _t, \kappa _2 +\theta _{\infty} ; \tilde{\lambda } ,\tilde{\mu } )$ to the line $L_{\infty} ^{*}$.
Let $(q_4, p_4)$ (resp. $(\tilde{q}_4, \tilde{p}_4)$) be the coordinate of $U_4 $ for the parameters $(\theta _0, \theta _1, \theta _t,\theta _{\infty}; \lambda ,\mu )$ (resp. $(\kappa _2 + \theta _0, \kappa _2 +\theta _1, \kappa _2 +\theta _t, \kappa _2 +\theta _{\infty} ; \tilde{\lambda } ,\tilde{\mu } )$ (see equations~(\ref{def:U01234i}), (\ref{eq:defrelSPI})).
Then we can express $\tilde{q}_4$ and $\tilde{p}_4$ by the variables $q_4$ and $p_4$. By setting $p_4 =0$, we have $\tilde{p}_4=0$ and $\tilde{q}_4 =q_4-\kappa_2 (t(\kappa _1 +\theta _t -1 )+\kappa _1 +\theta _1 -1 )$.
Hence the line $L_{\infty}^*$ with the parameter $(\theta _0, \theta _1, \theta _t,\theta _{\infty}; \lambda ,\mu )$ is mapped to the line $L_{\infty}^*$ with the parameter $(\kappa _2 + \theta _0, \kappa _2 +\theta _1, \kappa _2 +\theta _t, \kappa _2 +\theta _{\infty} ; \tilde{\lambda } ,\tilde{\mu } )$.
It follows from Proposition \ref{prop:Nov} that if $y_1(z)$ satisf\/ies equation~(\ref{eq:HeunP6infv2}) then the function $\tilde{y}(z)$ def\/ined by Proposition \ref{prop:Nov} satisf\/ies Heun's equation with the parameters $\gamma =1-\theta _0 -\kappa _2 $, $\delta = 1-\theta _1-\kappa _2$, $\epsilon = 1-\theta _t-\kappa _2$, $\alpha = 1- \kappa _2$, $\beta =1+ \theta _{\infty }$, $q=-q_4 -(1+t) \kappa _2 + (t(\kappa _1 +\theta _t)+\kappa _1+\theta _1))$, and the integral representation reproduces Theorem~\ref{thm:Heunint} by setting $\eta =\alpha = 1- \kappa _2$.
Therefore we have the following theorem:
\begin{theorem} \label{thm:k2Heun}
Let $X= L_0$, $L_1$, $L_t $ or $L_{\infty}^{*}$.
By the map $s_0s_3s_4s_2s_4s_3s_0: \: (\theta _0, \theta _1, \theta _t,\theta _{\infty}; \lambda ,\mu ) \mapsto (\kappa _2 + \theta _0, \kappa _2 +\theta _1, \kappa _2 +\theta _t, \kappa _2 +\theta _{\infty} ; \tilde{\lambda } ,\tilde{\mu } )$ where $\tilde{\lambda }$ and $\tilde{\mu } $ are defined in equation~\eqref{eq:tlamtmu}, the line~$X$ in the space of initial conditions with the parameter $(\theta _0, \theta _1, \theta _t,\theta _{\infty}; \lambda ,\mu )$ is mapped to the line~$X$ in the space of initial conditions with the parameter $(\kappa _2 + \theta _0, \kappa _2 +\theta _1, \kappa _2 +\theta _t, \kappa _2 +\theta _{\infty} ; \tilde{\lambda } ,\tilde{\mu } )$ where $\tilde{\lambda }$ and $\tilde{\mu } $ are defined in equation~\eqref{eq:tlamtmu},
and the integral transformation in Proposition~{\rm \ref{prop:Nov}} determined by the middle convolution reproduces the integral transformation on Heun's equations in Theorem~{\rm \ref{thm:Heunint}} by the restriction to the line~$X$.
\end{theorem}
Note that if $X= L^*_0$, $L^*_1$, $L^*_t $ or $L_{\infty}$ then the image of the line $X$ by the map $s_0s_3s_4s_2s_4s_3s_0 $ may not included in $X$ with the parameter $(\kappa _2 + \theta _0, \kappa _2 +\theta _1, \kappa _2 +\theta _t, \kappa _2 +\theta _{\infty} ; \tilde{\lambda } ,\tilde{\mu } )$.

We consider restrictions of the middle convolution for the case $\nu =\kappa _1$ (see Proposition \ref{prop:zinterepk1}) to lines in the space of initial conditions.
We discuss the restriction of the map $(\theta _0, \theta _1, \theta _t,\theta _{\infty}; \lambda ,\mu )$ $\mapsto (\kappa _1 + \theta _0, \kappa _1 +\theta _1, \kappa _1 +\theta _t, -\kappa _2 ; \lambda +\kappa _1 /\mu ,\mu  )$ to the line $L_0 ^{*}$.
Let $(q_1, p_1)$ (resp. $(\tilde{q}_1, \tilde{p}_1)$) be the coordinate of $U_ 1$ for the parameters $(\theta _0, \theta _1, \theta _t,\theta _{\infty}; \lambda ,\mu )$ (resp. $(\kappa _1 + \theta _0, \kappa _1 +\theta _1, \kappa _1 +\theta _t, -\kappa _2 ; \lambda +\kappa _1 /\mu ,\mu ))$.
Then we can express $\tilde{q}_1$ and $\tilde{p}_1$ by the variables $q_1$ and $p_1$, and by setting $p_1 =0$ we have $\tilde{p}_1=0$ and $\tilde{q}_1 = q_1 $.
Let $y_1 (z)$ be a solution to the Fuchsian equation $D_{y_1}(\theta _0, \theta _1, \theta _t,\theta _{\infty}; \lambda ,\mu )$ for the case $p_1=0 $ and set $v _1(z)=z^{-1} y_1(z)$.
Then $v_1 (z)$ satisf\/ies Heun's equation written as equation~(\ref{eq:HeunP60v2}).
On the case $p_1 =0$, the integral representation (equation~(\ref{eq:yt1zintrepk1})) is written as
\begin{gather*}
 \tilde{y}_1 (z)= \int _{[\gamma _{z} ,\gamma _p]} \left\{ \kappa _1 y_1(w) + w\frac{dy_1(w)}{dw}  \right\} \frac{(z-w)^{\kappa _1 }}{w} dw, \qquad (p\in \{0,1,t,\infty \}).
\end{gather*}
We set $\tilde{v} _1 (z)=z^{-1} \tilde{y} _1(z)$.
By integration by parts we have
\begin{gather}
  \tilde{v}_1 (z)= \frac{1}{z} \int _{[\gamma _{z} ,\gamma _p]} \left\{ \kappa _1 v_1(w) (z-w)^{\kappa _1 } + \frac{d(w v_1(w))}{dw} (z-w)^{\kappa _1 }   \right\} dw  \nonumber\\
\phantom{\tilde{v}_1 (z)}{}  = \frac{1}{z} \int _{[\gamma _{z} ,\gamma _p]} \big\{ \kappa _1 (z-w) v_1 (w) (z-w)^{\kappa _1 -1} + \kappa _1 w v_1 (w)  (z-w)^{\kappa _1 -1} \big\} dw \nonumber \\
\phantom{\tilde{v}_1 (z)}{}  = \kappa _1 \int _{[\gamma _{z} ,\gamma _p]}  v_1 (w) (z-w)^{\kappa _1 -1}  dw .\label{eq:intbyparts}
\end{gather}
On the other hand, it follows from Proposition \ref{prop:zinterepk1} and equation~(\ref{eq:HeunP60v2}) that $\tilde{v} _1 (z)$ satisf\/ies Heun's equation with the parameters $\gamma =2-\theta _0 -\kappa _1 $, $\delta = 1-\theta _1-\kappa _1$, $\epsilon = 1-\theta _t-\kappa _1$, $\alpha =1-\kappa _1$, $\beta =2-\theta _{\infty} $, $q=- tq_1 +(\kappa _1+\theta _0-1) \{ t(\kappa _1+\theta _1 -1) +\kappa _1+\theta _t-1 \}$.
Hence equation~(\ref{eq:intbyparts}) reproduces Theorem \ref{thm:Heunint} by setting $\eta = \alpha =1-\kappa _1$,
We can also obtain similar results for $L_1 ^*$, $L_t^*$ and $L_{\infty} ^*$.
Therefore we have the following theorem:
\begin{theorem} \label{thm:k1Heun}
Let $X= L_0^{*}$, $L_1^{*}$, $L_t ^{*}$ or $L_{\infty}^{*}$.
By the map $s_2: \: (\theta _0, \theta _1, \theta _t,\theta _{\infty}; \lambda ,\mu ) \mapsto (\kappa _1 + \theta _0,$ $\kappa _1 +\theta _1, \kappa _1 +\theta _t, -\kappa _2 ; \lambda +\kappa _1 /\mu ,\mu  )$, the line $X$ in the space of initial conditions with the parameter $(\theta _0, \theta _1, \theta _t,\theta _{\infty}; \lambda ,\mu )$ is mapped to the line $X$ in the space of initial conditions with the parameter $(\kappa _1 + \theta _0, \kappa _1 +\theta _1, \kappa _1 +\theta _t, -\kappa _2 ; \lambda +\kappa _1 /\mu ,\mu  )$,
and the integral transformation in Proposition~{\rm \ref{prop:zinterepk1}} determined by the middle convolution reproduces the integral transformation on Heun's equations in Theorem~{\rm \ref{thm:Heunint}} by the restriction to the line~$X$.
\end{theorem}

Note that if $X= L_0$, $L_1$, $L_t $ or $L_{\infty}$ then the image of the line $X$ by the map $s_2$ may not included in $X$ with the parameter $(\kappa _1 + \theta _0, \kappa _1 +\theta _1, \kappa _1 +\theta _t, -\kappa _2 ; \lambda +\kappa _1 /\mu ,\mu  )$.

\section{Middle convolution for the case that the parameter is integer} \label{sec:integer}

On the case $\kappa _2\in \Zint$, the function in equation~(\ref{eq:inttransDy1}) containing the Pochhammer contour may be vanished, and we propose other expressions of solutions to Fuchsian equation $D_{y_1}(\kappa _2 + \theta _0, \kappa _2 +\theta _1,$ $\kappa _2 +\theta _t, \kappa _2 +\theta _{\infty} ; \tilde{\lambda } ,\tilde{\mu } )$ in use of solutions to $D_{y_1}(\theta _0, \theta _1, \theta _t,\theta _{\infty}; \lambda ,\mu )$.

We have the following proposition for the case $\kappa _2\in \Zint _{<0}$,
\begin{proposition} \label{prop:diff}
$(i)$ Let $A_0$, $A_1$, $A_t$ be matrices in $\Cplx ^{2\times 2}$, and let $B^{(\nu )}_0, B^{(\nu )}_1, B^{(\nu )}_t \in \Cplx ^{6\times 6}$ be the matrices defined in equation~\eqref{eq:Bdef} for $\nu \in \Cplx$.
Assume that $\nu \in \Zint _{<0}$ and $Y= {}^t (y_1(z), y_2(z))$ is a~solution to the system of differential equations
\begin{gather*}
\frac{dY}{dz}=\left( \frac{A_0}{z}+\frac{A_1}{z-1}+\frac{A_t}{z-t} \right) Y.
%\label{eq:dYdzAzY00}
\end{gather*}
Write $\nu =-1-n$ $(n\in \Zint _{\geq 0})$.
Then the function
\begin{gather*}
U = \left(
\begin{array}{l}
(d/dz) ^{n} (z^{-1}y_1(z)) \\
(d/dz) ^{n} (z^{-1}y_2(z)) \\
(d/dz) ^{n} ((z-1)^{-1}y_1(z)) \\
(d/dz) ^{n} ((z-1)^{-1}y_2(z)) \\
(d/dz) ^{n} ((z-t)^{-1}y_1(z)) \\
(d/dz) ^{n} ((z-t)^{-1}y_2(z))
\end{array}
\right) , %\label{eq:diffrepU}
\end{gather*}
satisfies the system of differential equations
\begin{gather}
\frac{dU}{dz}=\left( \frac{B^{(-1-n )}_0}{z}+\frac{B^{(-1-n )}_1}{z-1}+\frac{B^{(-1-n )}_t}{z-t} \right) U.
\label{eq:dYdzBzU00}
\end{gather}

$(ii)$ If $\kappa _2 \in \Zint _{< 0}$ and $Y={}^t ( y_1(z), y_2(z))$ is a solution to the Fuchsian system $D_Y(\theta _0, \theta _1, \theta _t,\theta _{\infty};$ $\lambda ,\mu ;k )$,
then the function $\tilde{Y} = {}^t (\tilde{y}_1 (z) , \tilde{y}_2 (z) )$ defined by
\begin{gather}
  \tilde{y}_1 (z) = \left( \frac{d}{dz} \right) ^{-\kappa _2} y_1(z) , \nonumber \\
  \tilde{y}_2 (z)= \frac{\kappa _2 \lambda (\lambda -1) (\lambda -t)}{k(\lambda -\tilde{\lambda })} \left[ \left\{ \left( \frac{d}{dz} \right) ^{-\kappa _2} y_1(z) -\mu  \left( \frac{d}{dz} \right) ^{-\kappa _2-1} y_1(z) \right\} \frac{1}{\lambda -z} \right.\nonumber\\
  \left. \phantom{\tilde{y}_2 (z)=}{}
  +\frac{\mu}{\kappa _1}\left( \frac{d}{dz} \right) ^{-\kappa _2} y_1(z) \right], \label{eq:yt1zdiffrep}
\end{gather}
satisfies the Fuchsian system $D_Y(\kappa _2 + \theta _0, \kappa _2 +\theta _1, \kappa _2 +\theta _t, \kappa _2 +\theta _{\infty} ; \tilde{\lambda } ,\tilde{\mu };k )$, where $\tilde{\lambda }$ and $\tilde{\mu }$ are defined in equation~\eqref{eq:tlamtmu}.

$(iii)$ If $\kappa _2 \in \Zint _{< 0}$ and $y_1(z)$ is a solution to $D_{y_1}(\theta _0, \theta _1, \theta _t,\theta _{\infty}; \lambda ,\mu )$,
then the function
\begin{gather*}
\tilde{y} (z)  = \left( \frac{d}{dz} \right) ^{-\kappa _2} y_1(z), %\label{eq:yt1zdiffrepy1}
\end{gather*}
satisfies $D_{y_1}(\kappa _2 + \theta _0, \kappa _2 +\theta _1, \kappa _2 +\theta _t, \kappa _2 +\theta _{\infty} ; \tilde{\lambda } ,\tilde{\mu } )$.
\end{proposition}

\begin{proof}
$(i)$ If $\nu= -1$, then it follows immediately that the function $U={}^t ( z^{-1}y_1(z),z^{-1}y_2(z)$, $(z-1)^{-1}y_1(z),(z-1)^{-1}y_2(z),  (z-t)^{-1}y_1(z) ,(z-t)^{-1}y_2(z))$  satisf\/ied equation~(\ref{eq:dYdzBzU00}) for $n=0$.
Assume now that the function $U={}^t ( u_1(z),  u_2(z),  u_3(z),  u_4(z),  u_5(z),  u_6(z))$ satisf\/ies equation~(\ref{eq:dYdzBzU00}).
Set $V=dU/dz$.
Since
\begin{gather*}
  \frac{B^{(-1-n )}_0}{z^2}+\frac{B^{(-1-n )}_1}{(z-1)^2}+\frac{B^{(-1-n )}_t}{(z-t)^2}=
\frac{1}{z} \frac{B^{(-1-n )}_0}{z}
\left(
\begin{array}{c}
u_1(z) \\
u_2(z) \\
0 \\
0 \\
0 \\
0
\end{array}
\right)
+ \frac{1}{z-1} \frac{B^{(-1-n )}_0}{z-1}
\left(
\begin{array}{c}
0 \\
0 \\
u_3(z) \\
u_4(z) \\
0 \\
0
\end{array}
\right) \\
 {} +  \frac{1}{z-t} \frac{B^{(-1-n )}_0}{z-t}
\left(
\begin{array}{c}
0 \\
0 \\
0 \\
0 \\
u_5(z) \\
u_6(z)
\end{array}
\right)
=\frac{1}{z}
\left(
\begin{array}{c}
u'_1(z) \\
u'_2(z) \\
0 \\
0 \\
0 \\
0
\end{array}
\right)
+ \frac{1}{z-1}
\left(
\begin{array}{c}
0 \\
0 \\
u'_3(z) \\
u'_4(z) \\
0 \\
0
\end{array}
\right)
+  \frac{1}{z-t}
\left(
\begin{array}{c}
0 \\
0 \\
0 \\
0 \\
u'_5(z) \\
u'_6(z)
\end{array}
\right) , \nonumber
\end{gather*}
we have
\begin{gather*}
\frac{dV}{dz}  = \frac{d}{dz} \left\{ \left( \frac{B^{(-1-n )}_0}{z}+\frac{B^{(-1-n )}_1}{z-1}+\frac{B^{(-1-n )}_t}{z-t} \right) U \right\} \\
  \phantom{\frac{dV}{dz}}{} =-\left( \frac{B^{(-1-n )}_0}{z^2}+\frac{B^{(-1-n )}_1}{(z-1)^2}+\frac{B^{(-1-n )}_t}{(z-t)^2} \right) U+\left( \frac{B^{(-1-n )}_0}{z}+\frac{B^{(-1-n )}_1}{z-1}+\frac{B^{(-1-n )}_t}{z-t} \right) V \nonumber \\
\phantom{\frac{dV}{dz}}{}  =\left( \frac{B^{(-2-n )}_0}{z}+\frac{B^{(-2-n )}_1}{z-1}+\frac{B^{(-2-n )}_t}{z-t} \right) V.\nonumber
\end{gather*}
Hence $(i)$ is proved inductively.

$(ii)$ Let $A_0$, $A_1$, $A_t$ be the matrices def\/ined by equation~(\ref{eq:A0A1AtP}) and set $\nu = \kappa _2$.
We def\/ine the matrix $S$ by equation~(\ref{eq:matS}) and set $\tilde{U}= S^{-1} U$.
Then  $\tilde{Y}= {}^t (\tilde{u}_1(z), \tilde{u}_2(z))$ satisf\/ies the Fuchsian dif\/ferential equation $D_Y(\tilde{\theta }_0, \tilde{\theta }_1, \tilde{\theta }_t, \tilde{\theta }_{\infty}; \tilde{\lambda },\tilde{\mu };\tilde{k} )$ where the parameters are determined by equation~(\ref{eq:tt0tt1ttt}), and $\tilde{u}_1(z)$, $\tilde{u}_2(z)$ are expressed as equation~(\ref{eq:y1tu16}).
By a straightforward calculation as obtaining equation~(\ref{eq:yt1zintrep}), we have equation~(\ref{eq:yt1zdiffrep}).

$(iii)$ follows from $(ii)$.
\end{proof}

Note that $(i)$ is valid for Fuchsian dif\/ferential systems of arbitrary size and arbitrary number of regular singularities.
We have a similar statement for the case $\nu =\kappa _1$ and $\kappa _1 \in \Zint_{< 0}$.
Namely, if $\kappa _1 \in \Zint _{< 0}$ and $Y={}^t ( y_1(z), y_2(z))$ is a solution to the Fuchsian dif\/ferential equation $D_Y(\theta _0, \theta _1, \theta _t,\theta _{\infty}; \lambda ,\mu ;k )$ (see equation~(\ref{eq:dYdzAzY})), then the function $\tilde{Y} = {}^t (\tilde{y}_1 (z) , \tilde{y}_2 (z) )$ def\/ined by
\begin{gather*}
  \tilde{y}_1 (z)= \left( \frac{d}{dz} \right) ^{-\kappa _1 -1} \left\{ \frac{1}{z-\lambda } \left( \kappa _1 y_1(z) + (z- \tilde{\lambda})\frac{dy_1(z)}{dz}  \right) \right\}  ,  \nonumber\\
  \tilde{y}_2 (z)= \left( \frac{d}{dz} \right) ^{-\kappa _1 } y_2(z) ,  %\label{eq:yt1zdiffrepk1}
\end{gather*}
satisf\/ies the Fuchsian system $D_Y(\kappa _1 + \theta _0, \kappa _1 +\theta _1, \kappa _1 +\theta _t, -\kappa _2 ; \lambda +\kappa _1 /\mu ,\mu ;k )$.

If $\kappa _2=0$ (resp.\ $\kappa _1=0$), then the Fuchsian system $D_Y(\kappa _2 + \theta _0, \kappa _2 +\theta _1, \kappa _2 +\theta _t, \kappa _2 +\theta _{\infty} ; \tilde{\lambda } ,\tilde{\mu };k )$ (resp.\ $D_Y(\kappa _1 + \theta _0, \kappa _1 +\theta _1, \kappa _1 +\theta _t, -\kappa _2 ; \lambda +\kappa _1 /\mu ,\mu ;k )$) coincides with $D_Y(\theta _0, \theta _1, \theta _t,\theta _{\infty}; \lambda ,\mu ;k )$, and the function $\tilde{Y}={}^t ( \tilde{y}_1(z), \tilde{y}_2(z))$ just corresponds to ${}^t (y_1(z), y_2(z))$.

On the case $\kappa _2\in \Zint _{>0}$, we have the following proposition:
\begin{proposition}
Let $p\in \{ 0,1,t,\infty \}$ and $C_p$ be the cycle starting from $w=z$, turning $w=p$ anti-clockwise and return to $w=z$.

$(i)$ If $\kappa _2 \in \Zint _{> 0}$ and $Y={}^t ( y_1(z), y_2(z))$ is a solution to the Fuchsian system $D_Y(\theta _0, \theta _1, \theta _t,\theta _{\infty};$ $\lambda ,\mu ;k )$, then the function $\tilde{Y} = {}^t (\tilde{y}_1 (z) , \tilde{y}_2 (z) )$ defined by
\begin{gather}
  \tilde{y}_1 (z) = \int _{C_p}  \frac{dy_1(w)}{dw}  (z-w)^{\kappa _2 } dw, \label{eq:yt1zintrepk2pos} \\
  \tilde{y}_2 (z)= \frac{\kappa _2 \lambda (\lambda -1) (\lambda -t)}{k(\lambda -\tilde{\lambda })} \int _{C_p} \left\{ \left( \frac{dy_1(w)}{dw} -\mu y_1(w)  \right) \frac{1}{\lambda -w} +\frac{\mu}{\kappa _1}\frac{dy_1(w)}{dw} \right\} (z-w)^{\kappa _2 } dw,  \nonumber
\end{gather}
satisfies the Fuchsian system $D_Y(\kappa _2 + \theta _0, \kappa _2 +\theta _1, \kappa _2 +\theta _t, \kappa _2 +\theta _{\infty} ; \tilde{\lambda } ,\tilde{\mu };k )$ for $p \in \{ 0,1,t,\infty \}$ where $\tilde{\lambda }$ and $\tilde{\mu }$ are defined in equation~\eqref{eq:tlamtmu}.

$(ii)$ If $\kappa _2 \in \Zint _{> 0}$ and $y_1(z)$ is a solution to $D_{y_1}(\theta _0, \theta _1, \theta _t,\theta _{\infty}; \lambda ,\mu )$,
 then the function
\begin{gather}
\tilde{y} (z)  = \int _{C_p } y_1(w) (z-w)^{\kappa _2 -1} dw, \label{eq:yt1zrepy1k2pos}
\end{gather}
satisfies $D_{y_1}(\kappa _2 + \theta _0, \kappa _2 +\theta _1, \kappa _2 +\theta _t, \kappa _2 +\theta _{\infty} ; \tilde{\lambda } ,\tilde{\mu } )$ for $p \in \{ 0,1,t,\infty \}$.
\end{proposition}
Note that the functions $ \tilde{y}_1 (z)$, $ \tilde{y}_2 (z)$, $ \tilde{y} (z)$ may not be polynomials although the integrands of equations~(\ref{eq:yt1zintrepk2pos}) and (\ref{eq:yt1zrepy1k2pos}) are polynomials in $z$.
\begin{proof}
Set
\begin{gather*}
  K_1 (w) = \frac{dy_1(w)}{dw} , \\
 K_2 (w) = \frac{\kappa _2 \lambda (\lambda -1) (\lambda -t)}{k(\lambda -\tilde{\lambda })} \left\{ \left( \frac{dy_1(w)}{dw} -\mu y_1(w)  \right) \frac{1}{\lambda -w} +\frac{\mu}{\kappa _1}\frac{dy_1(w)}{dw} \right\} . \nonumber
\end{gather*}
It follows from Theorem \ref{thm:zinterep} that the function $Y(z)= \left( \! \begin{array}{c} \tilde{y}_1(z) \\ \tilde{y}_2(z) \end{array} \! \right)$ def\/ined by
\begin{gather*}
 \tilde{y}_i (z)  = \int _{\gamma _{z} \gamma _p \gamma _{z} ^{-1} \gamma _p ^{-1}} K_i(w) (z-w)^{\kappa _2 } dw \qquad \\
\phantom{\tilde{y}_i (z)}{}   = (1- e^{2\pi \sqrt{-1}\kappa _2}) \!\int _{\gamma _p} \! K_i (w) (z\!-\!w)^{\kappa _2 } dw + \! \int _{\gamma _z} \! (K^{\gamma _p} _i(w)\! - \! K_i(w))(z-w)^{\kappa _2 } dw   \quad (i=1,2)\!
\end{gather*}
is a solution to $D_Y(\kappa _2 + \theta _0, \kappa _2 +\theta _1, \kappa _2 +\theta _t, \kappa _2 +\theta _{\infty} ; \tilde{\lambda } ,\tilde{\mu };k )$ for $p \in \{0,1,t\}$, where $K^{\gamma _p} _i(w)$ is the function analytically continued along the cycle $\gamma _p$.

If $\kappa _2 >-1$, then the integrals $\int _{\gamma _z} ( K^{\gamma _p} _i(w) - K_i(w))(z-w)^{\kappa _2 } dw$ tends to zero as the base point $o$ of the integral tends to $z$,
 and it follows that the function $\left( \! \begin{array}{c} \int_{C_p } K_1(w) (z-w)^{\kappa _2 } dw \\  \int_{C_p } K_2(w) (z-w)^{\kappa _2 } dw \end{array} \! \right)$ is a solution to $D_Y(\kappa _2 + \theta _0, \kappa _2 +\theta _1, \kappa _2 +\theta _t, \kappa _2 +\theta _{\infty} ; \tilde{\lambda } ,\tilde{\mu };k )$ for $\kappa _2 \in \Rea _{>-1} \setminus \Zint _{>-1}$.
By considering the limit $\kappa _2 \rightarrow n$, $n\in \Zint_{>0}$, we obtain (i) for $p \in \{ 0,1,t\}$.
The case $ p=\infty $ follows from $C_{\infty} =C_t^{-1}C_1^{-1}C_0^{-1}$.

By integration by parts as equation~(\ref{eq:intbyparts0}) we obtain $(ii)$.
\end{proof}

We have similar proposition for the case $\nu =\kappa _1$ and $\kappa _1 \in \Zint _{>0}$.
On middle convolution $mc_\nu$ for Fuchsian dif\/ferential systems of arbitrary size and arbitrary number of regular singularities whose parameter $\nu $ is positive integer, the contour $[\gamma _z, \gamma _p]$ can be replaced by $C_p$.

We can reformulate Theorem~\ref{thm:k2Heun} (resp.\ Theorem~\ref{thm:k1Heun}) for the case $k_2 \in \Zint _{<0}$ (resp.\ $k_1 \in \Zint _{<0}$) by changing the integral to successive dif\/ferential and for the case $k_2 \in \Zint _{>0}$ (resp.\ $k_1 \in \Zint _{>0}$) by changing the contour $[\gamma _z ,\gamma _p] $ to $C_p$.
The corresponding setting for Heun's equation is described as follows:
\begin{proposition} \label{prop:Heuninteger}
Let $v(w)$ be a solution to Heun's equation written as equation~\eqref{Heun01}  with the parameters in equation~\eqref{eq:mualbe}.

$(i)$ If $\eta \in \Zint _{>1}$, then the function
\begin{gather*}
  y(z)= (d/dz)^{\eta -1} v(z)
\end{gather*}
is a solution to Heun's equation \eqref{Heun02}.

$(ii)$ If $\eta \in \Zint _{<1}$, then the function
\begin{gather*}
  y(z)=\int _{C_p} v(w) (z-w)^{-\eta } dw
\end{gather*}
is a solution to Heun's equation \eqref{Heun02}  for $p \in \{ 0,1,t,\infty \}$.
\end{proposition}

The generalized Darboux transformation (Crum--Darboux transformation) for elliptical representation of Heun's equation was introduced in~\cite{Tak5}, and we can show that Proposition~\ref{prop:Heuninteger}~$(i)$ gives another description of the generalized Darboux transformation.
Hence the integral transformation given by Theorem~\ref{thm:Heunint} can be regarded as a generalization of the generalized Darboux transformation to non-integer cases.
Khare and Sukhatme~\cite{KS0} conjectured a duality of quasi-exactly solvable (QES) eigenvalues for elliptical representation of Heun's equation.
By rewriting parameters of the duality to Heun's equation on the Riemann sphere, we obtain a correspondence on the parameters $\alpha$, $\beta$, $\gamma$, $\delta$, $\epsilon$, $q $ and $\alpha '$, $\beta '$, $\gamma '$, $\delta '$, $\epsilon '$, $q' $ on the integral transformation of Heun's equation in Theorem~\ref{thm:Heunint}.
We will report further from a viewpoint of monodromy in a~separated paper.

\appendix

\section{Appendix}

We investigate the realization of the Fuchsian system (equation~(\ref{eq:dYdzAzY1})) for the cases $\lambda = 0,1,t,\infty$ in the setting of Section~\ref{sec:linP6} and observe relationships with the lines $L_0$, $L^*_0$, $L_1$, $L^*_1$, $L_t$, $L^*_t$, $L_{\infty}$, $L^*_{\infty}$ (see equation~(\ref{eq:linesSIC})) in the space of initial conditions $E(t)$.

We consider the case $\lambda = 0$, i.e., the case $a_{12}^{(0)}= 0$, $a_{12}^{(1)}\neq 0$, $a_{12}^{(t)}\neq 0$, $(t+1)a_{12}^{(0)}+ t a_{12}^{(1)}+ a_{12}^{(t)}\neq  0 $ (see equation~(\ref{eq:Aiajk})).
Since $a_{12}^{(0)}= 0$ and the eigenvalues of $A_0$ are $\theta _0 $ and $0$, the matrix $A_0$ is written as
\begin{gather*}
A_0= \left(
\begin{array}{ll}
\theta _0  & 0 \\
v  & 0
\end{array}
\right)  \qquad
\mbox{or} \qquad
A_0= \left(
\begin{array}{ll}
0  & 0 \\
v  & \theta _0
\end{array}
\right) , %\label{eq:A00}
\end{gather*}
and it follows from $a_{12}^{(1)}\neq 0$, $a_{12}^{(t)}\neq 0 $ that the matrices $A_1$, $A_t$ may be expressed as equation~(\ref{eq:A0A1AtP}).
We determine $w_1$, $w_t$ so as to satisfy $a_{12} (z) =-w_1 (z-1)-w_t/(z-t) =k/(z(z-1))$.
Then we have
\begin{gather}
w_1=k/(t-1), \qquad w_t = -k/(t-1). \label{eq:w0}
\end{gather}

On the case
\begin{gather}
A_0= \left(
\begin{array}{ll}
0  & 0 \\
v  & \theta _0
\end{array}
\right) , \qquad
A_1= \left(
\begin{array}{ll}
u_1+\theta _1  & -w_1 \\
u_1(u_1+\theta _1)/w_1 & -u_1
\end{array}
\right) , \nonumber\\
A_t= \left(
\begin{array}{ll}
u_t+\theta _t  & -w_t \\
u_t(u_t+\theta _t)/w_t  & -u_t
\end{array}
\right), \label{eq:A001}
\end{gather}
we determine $u_1$, $u_t$ so as to satisfy equation~(\ref{def:Ainf}), namely
\begin{gather*}
  v +u_1(u_1+\theta _1)/ w_1 +u_t(u_t+\theta _t)/ w_t =0, \\
   u_1+\theta _1 + u_t+\theta _t =-\kappa _1, \quad  \theta _0 -u_1- u_t =-\kappa _2 . \nonumber
\end{gather*}
We have
\begin{gather}
  u_1=- \theta _1+ \frac{1}{\theta _{\infty }-\theta _0} \left( \frac{k v}{t-1}- \kappa _1 (\kappa _1 +\theta _t )\right)
 ,  \nonumber\\
  u_t= - \theta _t -\frac{1}{\theta _{\infty }-\theta _0} \left( \frac{k v}{t-1}+\kappa _1 (\kappa _1+\theta _1 )\right) .\label{eq:u0-1}
\end{gather}
Hence we have one-parameter realization of equation~(\ref{eq:A001}) with the prescribed condition.
We discuss relationship with the Fuchsian system on the line $L_0$.
For this purpose, we recall matri\-ces~$A_0$, $A_1$, $A_t$ determined by equations~(\ref{eq:A0A1AtP}), (\ref{eq:wgen}), (\ref{eq:ugen}) and restrict them to $\lambda =0$.
Then all elements in $A_0$, $A_1$ and $A_t$ are well-def\/ined and we have
\begin{gather*}
A_0 \left | _{\lambda =0 }\right. = \left(
\begin{array}{ll}
0  & 0 \\
\theta _ 0 \{ t(\theta _0-\theta _{\infty }) \mu+\kappa _ 1 (\kappa _1+\theta _1+t(\kappa _1+\theta _t))\} /(k\theta _{\infty }) & \theta _0
\end{array}
\right) .  % \label{eq:A0011}
\end{gather*}
In fact the matrices $A_0$, $A_1$, $A_t$ restricted to the line $L_0$ coincide with the ones determined by equations~(\ref{eq:A001}), (\ref{eq:w0}), (\ref{eq:u0-1}) and $v= \theta _ 0 \{ t(\theta _0-\theta _{\infty }) \mu+\kappa _ 1 (\kappa _1+\theta _1+t(\kappa _1+\theta _t))\} /(k\theta _{\infty })$.
Note that the second-order dif\/ferential equation for the function $y_1$ on the case of the matrices in equations~(\ref{eq:A001}), (\ref{eq:w0}), (\ref{eq:u0-1}) is obtained as equation~(\ref{eq:HeunP60}) by substituting
$\mu=(k\theta _{\infty } v-\kappa _ 1 \theta _0(\kappa _1+\theta _1+t(\kappa _1+\theta _t)))/(t\theta _0(\theta _0-\theta _{\infty })) $.

On the case
\begin{gather}
A_0= \left(
\begin{array}{ll}
\theta _0  & 0 \\
v  & 0
\end{array}
\right) , \qquad
A_1= \left(
\begin{array}{ll}
u_1+\theta _1  & -w_1 \\
u_1(u_1+\theta _1)/w_1 & -u_1
\end{array}
\right) , \nonumber\\
A_t= \left(
\begin{array}{ll}
u_t+\theta _t  & -w_t \\
u_t(u_t+\theta _t)/w_t  & -u_t
\end{array}
\right) , \label{eq:A002}
\end{gather}
$u_1$,  $u_t$ are determined as
\begin{gather}
  u_1=\frac{1}{\theta _{\infty }+ \theta _0} \left( \frac{kv}{t-1}-\kappa _2(\kappa _2+\theta _t) \right)  , \; u_t= \frac{-1}{\theta _{\infty }+ \theta _0} \left( \frac{kv}{t-1}+\kappa _2(\kappa _2+\theta _1) \right) , \label{eq:u0-2}
\end{gather}
to satisfy equation~(\ref{def:Ainf}).
To realize the Fuchsian system on the line $L^*_0$, we recall matrices $A_0$, $A_1$, $A_t$ determined by equations~(\ref{eq:A0A1AtP}), (\ref{eq:wgen}), (\ref{eq:ugen}), transform $(\lambda ,\mu) (=(q_0,p_0) )$ to $(q_1,p_1)$ by equation~(\ref{eq:defrelSPI}) and restrict matrix elements to $q_1 =0$.
Then the matrices $A_0$, $A_1$ and $A_t$ are determined as equations~(\ref{eq:A002}), (\ref{eq:w0}), (\ref{eq:u0-2}), where
\begin{gather}
v= \{ -t(\theta _0+\theta _{\infty })q_1+\theta _0 (\kappa _1+\theta _0)((\kappa _2+\theta _t)+t(\kappa _2+\theta _1))\} /(k\theta _{\infty })
 .  \label{eq:A0012}
\end{gather}
Note that the second-order dif\/ferential equation for the function $\tilde{y}_1 =z^{-1}y_1$ on the case of the matrices in equations~(\ref{eq:A002}), (\ref{eq:w0}), (\ref{eq:u0-2}) is obtained as equation~(\ref{eq:HeunP60v2}) by substituting equation~(\ref{eq:A0012}).

We consider the case $\lambda = 1$, i.e., the case $a_{12}^{(1)}= 0 , a_{12}^{(0)}\neq 0 , a_{12}^{(t)}\neq 0 , (t+1)a_{12}^{(0)}+ t a_{12}^{(1)}+ a_{12}^{(t)}\neq  0 $.
Since $a_{12}^{(1)}= 0$ and the eigenvalues of $A_1$ are $\theta _1$ and $0$, the matrix $A_1$ is written as
\begin{gather*}
A_1= \left(
\begin{array}{ll}
\theta _1  & 0 \\
v  & 0
\end{array}
\right)  \qquad
\mbox{or} \qquad
A_1= \left(
\begin{array}{ll}
0  & 0 \\
v  & \theta _1
\end{array}
\right) , %\label{eq:A10}
\end{gather*}
and the matrices $A_0$, $A_t$ may be expressed as equation~(\ref{eq:A0A1AtP}).
To satisfy $a_{12} (z) =-w_0/z-w_t/(z-t) =k/(z(z-t))$, we have
\begin{gather}
w_0=k/t, \qquad w_t = -k/t. \label{eq:w1}
\end{gather}
On the case
\begin{gather}
A_1= \left(
\begin{array}{ll}
0  & 0 \\
v  & \theta _1
\end{array}
\right) ,  \qquad
A_0= \left(
\begin{array}{ll}
u_0+\theta _0  & -w_0 \\
u_0(u_0+\theta _0)/w_0 & -u_0
\end{array}
\right) , \nonumber\\
A_t= \left(
\begin{array}{ll}
u_t+\theta _t  & -w_t \\
u_t(u_t+\theta _t)/w_t  & -u_t
\end{array}
\right) ,
\label{eq:A101}
\end{gather}
$u_0$, $u_t$ are determined as
\begin{gather}
  u_0=- \theta _0+ \frac{1}{\theta _{\infty }-\theta _1} \left( \frac{k v}{t}- \kappa _1 (\kappa _1 +\theta _t )\right)
 , \nonumber\\
  u_t= - \theta _t - \frac{1}{\theta _{\infty }-\theta _1} \left( \frac{k v}{t}+\kappa _1 (\kappa _1+\theta _0 )\right) ,\label{eq:u1-1}
\end{gather}
to satisfy equation~(\ref{def:Ainf}).
To realize the Fuchsian system on the line $L_1$, we recall matrices $A_0$, $A_1$, $A_t$ determined by equations~(\ref{eq:A0A1AtP}), (\ref{eq:wgen}), (\ref{eq:ugen}) and restrict matrix elements to $\lambda =1$.
Then the matrices $A_0$, $A_1$ and $A_t$ are determined as equations~(\ref{eq:A101}), (\ref{eq:w1}), (\ref{eq:u1-1}), where
\begin{gather}
v= \theta _ 1 \{ (1-t)(\theta _1-\theta _{\infty }) \mu-\kappa _ 1 (\kappa _1+\theta _0+(1-t)(\kappa _1+\theta _t))\}/(k\theta _{\infty }).  \label{eq:A0021}
\end{gather}
Note that the second-order dif\/ferential equation for the function $y_1$ on the case of the matrices in equations~(\ref{eq:A101}), (\ref{eq:w1}), (\ref{eq:u1-1}) is obtained as equation~(\ref{eq:HeunP61}) by substituting equation~(\ref{eq:A0021}).

On the case
\begin{gather}
A_1= \left(
\begin{array}{ll}
\theta _1  & 0 \\
v  & 0
\end{array}
\right) , \qquad
A_0= \left(
\begin{array}{ll}
u_0+\theta _0  & -w_0 \\
u_0(u_0+\theta _0)/w_0 & -u_0
\end{array}
\right) , \nonumber\\
A_t= \left(
\begin{array}{ll}
u_t+\theta _t  & -w_t \\
u_t(u_t+\theta _t)/w_t  & -u_t
\end{array}
\right) ,
\label{eq:A102}
\end{gather}
$u_0$,  $u_t$ are determined as
\begin{gather}
  u_0=\frac{1}{\theta _{\infty }+ \theta _1} \left( \frac{kv}{t}-\kappa _2(\kappa _2+\theta _t) \right)  , \qquad u_t= \frac{-1}{\theta _{\infty }+ \theta _1} \left( \frac{kv}{t}+\kappa _2(\kappa _2+\theta _0) \right) , \label{eq:u1-2}
\end{gather}
to satisfy equation~(\ref{def:Ainf}).
To realize the Fuchsian system on the line $L^*_1$, we recall matrices $A_0$, $A_1$, $A_t$ determined by equations~(\ref{eq:A0A1AtP}), (\ref{eq:wgen}), (\ref{eq:ugen}), transform $(\lambda ,\mu) (=(q_0,p_0) )$ to $(q_2,p_2)$ by equation~(\ref{eq:defrelSPI}) and restrict matrix elements to $q_2 =0$.
Then the matrices $A_0$, $A_1$ and $A_t$ are determined as equations~(\ref{eq:A102}), (\ref{eq:w1}), (\ref{eq:u1-2}), where
\begin{gather}
v= \{ (t-1) (\theta _1+\theta _{\infty })q_2-\theta _1 (\kappa _1+\theta _1)((\kappa _2+\theta _t)+(1-t)(\kappa _2+\theta _0))\}/(k\theta _{\infty }) .  \label{eq:A0022}
\end{gather}
Note that the second-order dif\/ferential equation for the function $\tilde{y}_1 =(z-1)^{-1}y_1$ on the case of the matrices in equations~(\ref{eq:A102}), (\ref{eq:w1}), (\ref{eq:u1-2}) is obtained as equation~(\ref{eq:HeunP61v2}) by substituting equation~(\ref{eq:A0022}).

We consider the case $\lambda = t$, i.e., the case $a_{12}^{(t)}= 0$, $a_{12}^{(0)}\neq 0$, $a_{12}^{(1)}\neq 0$, $(t+1)a_{12}^{(0)}+ t a_{12}^{(1)}+ a_{12}^{(t)}\neq  0 $.
Then the matrix $A_t$ is written as
\begin{gather*}
A_t= \left(
\begin{array}{ll}
\theta _t  & 0 \\
v  & 0
\end{array}
\right)  \qquad
\mbox{or} \qquad
A_t= \left(
\begin{array}{ll}
0  & 0 \\
v  & \theta _t
\end{array}
\right) , %\label{eq:At0}
\end{gather*}
and the matrices $A_0$, $A_1$ may be expressed as equation~(\ref{eq:A0A1AtP}).
To satisfy $a_{12} (z) =-w_0/z-w_1/(z-1) =k/(z(z-1))$, we have
\begin{gather}
w_0=k, \qquad w_1 = -k. \label{eq:wt}
\end{gather}
On the case
\begin{gather}
A_t= \left(
\begin{array}{ll}
0  & 0 \\
v  & \theta _t
\end{array}
\right) ,  \qquad
A_0= \left(
\begin{array}{ll}
u_0+\theta _0  & -w_0 \\
u_0(u_0+\theta _0)/w_0 & -u_0
\end{array}
\right) , \nonumber\\
A_1= \left(
\begin{array}{ll}
u_1+\theta _1  & -w_1 \\
u_1(u_1+\theta _1)/w_1  & -u_1
\end{array}
\right) ,
\label{eq:At01}
\end{gather}
$u_0$, $u_1$ are determined as
\begin{gather}
  u_0=- \theta _0+ \frac{1}{\theta _{\infty }-\theta _t} \left( k v- \kappa _1 (\kappa _1 +\theta _1 )\right)
  , \nonumber\\
   u_1= - \theta _1 - \frac{1}{\theta _{\infty }-\theta _t} \left( k v + \kappa _1 (\kappa _1+\theta _0 )\right) ,
\label{eq:ut-1}
%\nonumber
\end{gather}
to satisfy equation~(\ref{def:Ainf}).
To realize the Fuchsian system on the line $L_t$, we recall matrices $A_0$, $A_1$, $A_t$ determined by equations~(\ref{eq:A0A1AtP}), (\ref{eq:wgen}), (\ref{eq:ugen}) and restrict matrix elements to $\lambda =t$.
Then the matrices $A_0$, $A_1$ and $A_t$ are determined as equations~(\ref{eq:At01}), (\ref{eq:wt}), (\ref{eq:ut-1}), where
\begin{gather}
v= \theta _ t \{ t(t-1)(\theta _t-\theta _{\infty }) \mu-\kappa _ 1 (t(\kappa _1+\theta _0)+(t-1)(\kappa _1+\theta _1))\} /(k\theta _{\infty }) .  \label{eq:A0031}
\end{gather}
Note that the second-order dif\/ferential equation for the function $y_1$ on the case of the matrices in equations~(\ref{eq:At01}), (\ref{eq:wt}), (\ref{eq:ut-1}) is obtained as equation~(\ref{eq:HeunP6t}) by substituting equation~(\ref{eq:A0031}).

On the case
\begin{gather}
A_t= \left(
\begin{array}{ll}
\theta _t  & 0 \\
v  & 0
\end{array}
\right) , \qquad
A_0= \left(
\begin{array}{ll}
u_0+\theta _0  & -w_0 \\
u_0(u_0+\theta _0)/w_0 & -u_0
\end{array}
\right) , \nonumber\\
A_1= \left(
\begin{array}{ll}
u_1+\theta _1  & -w_1 \\
u_1(u_1+\theta _1)/w_1  & -u_1
\end{array}
\right) ,
\label{eq:At02}
\end{gather}
$u_0$,  $u_t$ are determined as
\begin{gather}
  u_0=\frac{1}{\theta _{\infty }+ \theta _t} \left( kv-\kappa _2(\kappa _2+\theta _1) \right)  , \qquad
   u_1= \frac{-1}{\theta _{\infty }+ \theta _t} \left( kv+\kappa _2(\kappa _2+\theta _0) \right) , \label{eq:ut-2}
\end{gather}
to satisfy equation~(\ref{def:Ainf}).
To realize the Fuchsian system on the line $L^*_t$, we recall matrices $A_0$, $A_1$, $A_t$ determined by equations~(\ref{eq:A0A1AtP}), (\ref{eq:wgen}), (\ref{eq:ugen}), transform $(\lambda ,\mu) (=(q_0,p_0) )$ to $(q_3,p_3)$ by equation~(\ref{eq:defrelSPI}) and restrict matrix elements to $q_3 =0$.
Then the matrices $A_0$, $A_1$ and $A_t$ are determined as equations~(\ref{eq:At02}), (\ref{eq:wt}), (\ref{eq:ut-2}), where
\begin{gather}
v= \{ t(1-t)(\theta _t +\theta _{\infty })q_3-\theta _t (\kappa _1+\theta _t)(t(\kappa _2+\theta _1)+(t-1)(\kappa _2+\theta _0)) \} /(k\theta _{\infty }).  \label{eq:A0032}
\end{gather}
Note that the second-order dif\/ferential equation for the function $\tilde{y}_1 =(z-1)^{-1}y_1$ on the case of the matrices in equations~(\ref{eq:At02}), (\ref{eq:wt}), (\ref{eq:ut-2}) is obtained as equation~(\ref{eq:HeunP6tv2}) by substituting equation~(\ref{eq:A0032}).

We consider the case $\lambda = \infty$, i.e., the case $a_{12}^{(0)}\neq 0$, $a_{12}^{(1)}\neq 0$, $a_{12}^{(t)}\neq 0 , (t+1)a_{12}^{(0)}+ t a_{12}^{(1)}+$ $a_{12}^{(t)}= 0 $.
We can set $A_0$, $A_1$, $A_t$ as equation~(\ref{eq:A0A1AtP}) and we determine $u_0$,  $u_1$,  $u_t$,  $w_0$,  $w_1$,  $w_t$ to satisfy equation~(\ref{def:Ainf}) and
\begin{gather*}
  a_{12} (z) =-\frac{w_0}{z}-\frac{w_1}{z-1}-\frac{w_t}{z-t} =\frac{k }{z(z-1)(z-t)}.
\end{gather*}
Then we have
\begin{gather}
w_0 = - k/t, \qquad w_1=k/(t-1), \qquad w_t = -k/(t(t-1)), \label{eq:wi}
\end{gather}
and other relations are written as
\begin{gather*}
  u_0(u_0+\theta _0)/ w_0 +u_1(u_1+\theta _1)/ w_1 +u_t(u_t+\theta _t)/ w_t =0, \\
  u_0+\theta _0 + u_1+\theta _1 + u_t+\theta _t =-\kappa _1, \qquad  -u_0-u_1- u_t =-\kappa _2 . \nonumber
\end{gather*}
We solve the equation for $u_0 $, $u_1 $, $u_t$ by adding one more relation $-( u_1+\theta _1 + t(u_t+\theta _t ))= l$.
We have
\begin{gather}
  u_0= -\theta _0 -\kappa _1 + \frac{\tilde{l}}{t \theta _{\infty}}  , \qquad
 u_1= -\theta _1 + \frac{\tilde{l} - \theta _{\infty } l}{(1-t)  \theta _{\infty}}  , \qquad
 u_t= -\theta _t + \frac{\tilde{l} - t\theta _{\infty } l }{t(t-1)  \theta _{\infty}}  , \nonumber \\
  \tilde{l} =  l^2+(\theta _1 +t\theta _t) l+ t\kappa _1(\kappa _1+\theta _0) .  \label{eq:ui}
\end{gather}
The second-order dif\/ferential equation for the function $y_1$ is written as
\begin{gather*}
 \frac{d^2y_1}{dz^2} + \left( \frac{1-\theta _0}{z}+\frac{1-\theta _1}{z-1}+\frac{1-\theta _t}{z-t} \right)  \frac{dy_1}{dz}+  \frac{\kappa _1 (\kappa _2 +2) z- q}{z(z-1)(z-t)} y_1=0,  \nonumber \\
  q=l (\theta _{\infty }-1) +\kappa _1 ( t(\kappa _2  +\theta _t +1)+\kappa _2 +\theta _1 +1).%\label{eq:HeunP6inf}
\end{gather*}
To realize the Fuchsian system on the line $L_{\infty }$, we recall matrices $A_0$, $A_1$, $A_t$ determined by equations~(\ref{eq:A0A1AtP}), (\ref{eq:wgen}), (\ref{eq:ugen}), transform $(\lambda ,\mu) (=(q_0,p_0) )$ to $(q_{\infty },p_{\infty })$ by equation~(\ref{eq:defrelSPI}), replace~$k$ by $-k q_{\infty } $ and restrict matrix elements to $q_{\infty } =0$.
Then the matrices $A_0$, $A_1$ and $A_t$ are determined as equations~(\ref{eq:A0A1AtP}), (\ref{eq:wi}), (\ref{eq:ui}), where
\begin{gather*}
l = p_{\infty } .
\end{gather*}

We have observed that Fuchsian systems on the lines $L_0$, $L^*_0$, $L_1$, $L^*_1$, $L_t$, $L^*_t$, $L_{\infty}$ are realized by Fuchsian systems for the case $\lambda = 0,1,t,\infty$.
Here the case $L^{*} _{\infty} $ is missing.
In fact this case does not simply correspond with the Fuchsian system as equation~(\ref{eq:dYdzAzY1}), because we have
\begin{gather*}
u_0= \frac{1}{t\theta _{\infty}}\frac{1}{p_4^2} +O\left(p_4^{-1}\right), \qquad  u_1= \frac{1}{(1-t)\theta _{\infty}}\frac{1}{p_4^2} +O\left(p_4^{-1}\right), \\ u_t= \frac{1}{t(t-1)\theta _{\infty}}\frac{1}{p_4^2} +O\left(p_4^{-1}\right),
\end{gather*}
as $p_4 \rightarrow 0 $ in the coordinate $(q_4, p_4)$ in equation~(\ref{eq:defrelSPI}),
although we can restrict the second-order dif\/ferential equation for $y_1$ to $p_4=0$ as equation~(\ref{eq:HeunP6infv2}).

\subsection*{Acknowledgements}
%{\bf Acknowledgments}
The author would like to thank Alexander Kazakov for sending the paper \cite{KS1} with valuable comments.
The author is supported by the Grant-in-Aid for Young Scientists (B) (No. 19740089) from the Japan Society for the Promotion of Science.

\pdfbookmark[1]{References}{ref}
\LastPageEnding

\end{document}